\documentclass[12pt]{article}
\usepackage{mathrsfs}
\usepackage{amsfonts}
\usepackage{amsmath,amssymb,amsfonts,amsthm,fancyhdr,mathrsfs}
\usepackage{appendix}
\usepackage{enumitem}
\usepackage{bbm}
\usepackage{verbatim}
\usepackage{color}
\usepackage{dsfont}
\usepackage[left=2.4cm,right=2.4cm,top=3cm,bottom=3cm]{geometry}

\usepackage{latexsym}
\usepackage{graphicx}
\usepackage{enumerate}
\usepackage[colorlinks=true,urlcolor=blue,
citecolor=red,linkcolor=blue,linktocpage,pdfpagelabels,
bookmarksnumbered,bookmarksopen]{hyperref}
\usepackage{cite}
\usepackage{indentfirst}
\usepackage{amsmath}
\allowdisplaybreaks[4]
\numberwithin{equation}{section} 

\numberwithin{equation}{section}
\newtheorem{thm}{Theorem}[section]

\newtheorem{pro}[thm]{Proposition}
\newtheorem{lem}[thm]{Lemma}

\newtheorem{re}{Remark}[section]
\newenvironment{pf}{{\noindent \it \bf Proof:}}{{\hfill$\Box$}\\}

\title{The Ground State of  a Cubic-quintic Nonlinear Schr\"{o}dinger Equation with Radial Potential in the Thomas-Fermi Limit}
\author{ \bf
	{\large Deke Li\thanks{School of Mathematics and Statistics, Lanzhou University, Lanzhou 730000, Gansu, People’s Republic of China. \emph{E-mail}: lidk2024@lzu.edu.cn.
	}},\quad
	{\large Qingxuan Wang\thanks{School of Mathematical Sciences, 
	Zhejiang Normal University, Jinhua 321004, Zhejiang, People’s Republic of China. \emph{E-mail}: wangqx@zjnu.edu.cn.} \thanks{Qingxuan Wang is the corresponding author.}}
}

\date{}
\begin{document}
\maketitle
\begin{abstract}
 We focus on  the ground state of  the   cubic-quintic nonlinear Schr\"{o}dinger energy functional 
\begin{gather*}
    \begin{aligned} 
		{E}(\varphi)=\frac{1}{2}\int_{\mathbb{R}^d}\left(|\nabla \varphi|^2+V(x)|\varphi|^2\right)\,dx
		\pm\frac{1}{4}\int_{\mathbb{R}^d}|\varphi|^4\,dx
		+\frac{1}{6}\int_{\mathbb{R}^d}|\varphi|^6\,dx, (d=1,2,3)
    \end{aligned}
\end{gather*}
 under the mass constraint $\int_{\mathbb{R}^d}|\varphi|^2\,dx=N$, where $N$ can be viewed as particle number, and $V(x)$ behaves like $C|x|^p (p\geq 2)$ as $|x|\rightarrow +\infty$, including the harmonic potential. When $N\rightarrow +\infty$, we show that up to a suitable scaling the ground state $\varphi_N$ would convergence strongly in some $L^q(\mathbb{R}^d)$ space to a Thomas-Fermi minimizer, this limit can be referred to as the \emph{Thomas-Fermi limit}. The limit Thomas-Fermi profile has compact support,  given by $u^{TF}(x)=\left[\mu^{TF}-C_0|x|^p\right]^{\frac{1}{4}}_{+}$, where $\mu^{TF}$ is a suitable Lagrange multiplier with exact value.  We find that, similar to the asymptotic analysis  in [J. Funct. Anal. 260 (2011), 2387-2406.] and [Arch. Ration. Mech. Anal. 217 (2015), 439-523.] for Gross-Pitaevskii energy in the Thomas-Fermi limit where a small parameter $\varepsilon$ tends to 0,  there also has a steep \emph{corner layer} near the boundary of compact support of $u^{TF}(x)$, in which the ground state has irregular behavior as $N\rightarrow +\infty$. Finally, we  establish a new energy method to obtain  the $L^\infty$-convergence rates of ground states $\varphi_N$ inside the corner layer and outside corner layer respectively, this method  may be applicable to other general nonlinearities.
\end{abstract}

%\vspace{0.2cm}



	\tableofcontents

\section{Introduction}
\subsection{The Problem}% and Mian Results}
This paper focuses on  the $N\to+\infty$ limiting behavior  of the cubic-quintic nonlinear Schr\"{o}dinger energy functional 
\begin{gather}\label{11}
	\begin{aligned} 
		{E}(\varPhi)=\frac{1}{2}\int_{\mathbb{R}^d}\left(|\nabla \varPhi|^2+V(x)|\varPhi|^2\right)\,dx
		+\frac{\kappa}{4}\int_{\mathbb{R}^d}|\varPhi|^4\,dx
		+\frac{1}{6}\int_{\mathbb{R}^d}|\varPhi|^6\,dx
	\end{aligned}
\end{gather}
minimized in
\begin{gather}\label{HN}
	\begin{aligned} 
		\mathscr{H}_N=\left\lbrace \varPhi\in H^1(\mathbb{R}^d):\ \int_{\mathbb{R}^d}V(x)|\varPhi|^2\,dx<\infty, \int_{\mathbb{R}^d}|\varPhi|^2\,dx=N\right\rbrace.
	\end{aligned}
\end{gather}
Where $d=1,2,3$,  the parameter $\kappa=\pm1$ ($\kappa> 0$ denotes defocusing interaction and $\kappa< 0$ denotes focusing interaction), $N>0$ denotes the particle number and the external potential $V(x)$  satisfies the following conditions:

\noindent\textbf{(V$_1$)}\quad $V\in C^1$ is increasing and non-negative radially symmetric;\\
\textbf{(V$_2$)}\quad There exists a  constant $C_0>0$ and $p\geq2$ such that
\begin{align}\label{188}
	\lim_{|x|\to\infty}\frac{\nabla V(x)\cdot x}{p|x|^p}=C_0;
\end{align}	
\textbf{(V$_3$)}\quad There exists constants $C_1\geq0$ and $0<\alpha<p$ such that
\begin{align}\label{HV}
	\lim_{|x|\to\infty}\frac{\nabla V(x)\cdot x-C_0p|x|^p}{|x|^\alpha}=C_1;
\end{align}
\textbf{(V$_4$)}\quad There exists a constant $C_2\geq0$ such that
\begin{align}\label{144}
	\lim_{|x|\to\infty}\frac{V(x)-C_0|x|^p}{|x|^\alpha}=C_2.
\end{align}
%Thus, it is common to refer to the above problem as the minimization of $E(\cdot)$ under the $L^2$-mass constraint.

There are many classes of external potentials $V$ satisfying the above conditions, such as harmonic external potential $|x|^2$ and more general external potentials of polynomial type $\Sigma_{i=1}^na_i|x|^{p_i}$ ($a_i>0,p_i\geq2$). In addition, in recent experiments, in which a laser beam is superimposed upon the magnetic trap holding the atoms, the trapping potential $V$ is of a different type \cite{Ryu07,Stock04}:
\begin{gather}\label{magnetic}
	\begin{aligned} 
		V(x)=|x|^2+ae^{-b|x|^2},\quad a,b>0.
	\end{aligned}
\end{gather}
By choosing $a,b$ accordingly, the magnetic trap satisfies \textbf{(V$_1$)}-\textbf{(V$_4$)}.

\subsection{Background and Motivation}
The cubic-quintic Schr\"{o}dinger equation associated with the energy functional \eqref{11} is given by
\begin{align}\label{1}
i\partial_t\Psi=-\Delta\Psi+V(x)\Psi+\kappa|\Psi|^2\Psi+|\Psi|^4\Psi,\quad(t,x)\in\mathbb{R}^+\times\mathbb{R}^d.
\end{align} 
This equation appears in numerous problems in physics, this can be used to describe nonlinear optics, the field theory, the mean-field theory of superconductivity, the motion of Bose-Einstein condensates, and Langmuir waves in plasma physics; one can see \cite{PREdesvataikov,PREchen,PREmih}  and the references therein.  The cubic term (also known as Kerr nonlinearity \cite{PREdesvataikov}) with a negative (positive) coefficient is a very natural case in physics. In particular, the incorporation of the defocusing quintic term is motivated by the stabilization of two-dimensional or three-dimensional vortex solitons \cite{PDmalomed}. Atoms/molecules do experience van der Waals attraction at intermediate densities. In fact, this attraction is the basis for the condensation of gases into liquids at low temperatures.  

The cubic-quintic  Schr\"{o}dinger equation is still getting a lot of attention in mathematics in the last decade. Recently, Martel \cite{IMMartel} studied the asymptotic stability of small standing solitary waves of the Schr\"odinger equation with a cubic-quintic nonlinearity in one space dimension, similar conclusions are presented in \cite{PMPMartel,CazenaveCMP,Ohta} and the references therein. Nguyen and Ricaud \cite{Nguyen-Ricaud24,Nguyen-Ricaud25} derived rigorously the cubic-quintic nonlinear Schr\"odinger semiclassical theory as the mean-field limit of the model and they investigated the behavior of the system in a double-limit. Killip,  Murphy, Visan and their collaborators \cite{Killip-17,Killip-18,Killip-21} studied the solitons, scattering and the initial-value problem with non-vanishing boundary conditions for the cubic-quintic nonlinear Schr\"odinger equation on $\mathbb{R}^3$, also refer to \cite{Murphy,Ardila23} and the references therein. 
The uniqueness, non-degeneracy of the positive solution for time-independent cubic-quintic nonlinear Schr\"odinger equation was shown in \cite{Carles-21,Lewin-20}. 
Soave \cite{Soave-jde,Soave-jfa} studied the normalized ground states for the nonlinear Schr\"odinger equation with combined nonlinearities including cubic-quintic nonlinearity.

A standing waves of \eqref{1} are solutions of the form
\begin{gather}\label{form}
	\begin{aligned} 
		\Psi(t,x)=e^{-i\mu t}\varPhi(x),
	\end{aligned}
\end{gather}
where $\mu\in\mathbb{R}$ and $\varPhi(x)\in H^1(\mathbb{R}^d)$ is a time-independent function. Usually, $\varPhi(x)$ is called a \emph{ground state} if it is a minimizer of following minimizing problem under the prescribed $L^2$ mass
\begin{gather}\label{m}
	\begin{aligned} 
		e(N):=\inf\left\lbrace E(\varPhi):\varPhi\in \mathscr{H}_N\right\rbrace,
	\end{aligned}
\end{gather}
where the energy functional $E(\varPhi)$ be given by \eqref{11}. The existence of ground states for this energy is stated as follow.

\begin{thm}\label{existence}
Suppose that $1\leq d\leq3$, $V(x)$ satisfy \textbf{(V$_1$)} and \textbf{(V$_4$)}. For any $N>0$, there exists at least one non-negative radially symmetric ground state for $e(N)$.
\end{thm}

The proof of this theorem  follows standard arguments and we include it here for completeness,  see  Appendix \ref{secA}.

\textbf{The first motivation} of this paper comes from the work by Lieb, Seiringer and Yngvason \cite{PRALieb,CMPLieb}, where they considered the following Gross-Pitaevskii energy functional
\begin{gather}\label{GPenergy}
	\begin{aligned} 
		E^{GP}(\varPhi)=\int_{\mathbb{R}^3}\left(|\nabla \varPhi|^2+V(x)|\varPhi|^2\right)\,dx
		+4\pi a\int_{\mathbb{R}^3}|\varPhi|^4\,dx,
	\end{aligned}
\end{gather}
where $a>0$ denotes the scattering length. For a given $N$ the corresponding Gross-Pitaevskii energy, denoted $\mathscr{E}^{GP}(N,a)$, is defined as the infimum of $E^{GP}(\varPhi)$ under the normalization condition
\begin{gather}\label{normal}
	\begin{aligned} 
		\int_{\mathbb{R}^3}|\varPhi|^2\,dx=N.
	\end{aligned}
\end{gather}
It has the simple scaling property
\begin{gather}\label{scaling}
	\begin{aligned} 
		\mathscr{E}^{GP}(N,a)=N\mathscr{E}^{GP}(1,Na).
	\end{aligned}
\end{gather}
What Eq. \eqref{scaling} shows is that the Gross-Pitaevskii energy functional \eqref{GPenergy} together with the normalization condition \eqref{normal} has one characteristic parameter, namely, $Na$. They consider the nontrivial aspects of Gross-Pitaevskii theory as $Na\to+\infty$.  It is sometimes referred to as the \emph{Thomas-Fermi limit} of Gross-Pitaevskii theory because the gradient term vanishes in this limit. They showed that the ground states $\varPhi_N$  of $\mathscr{E}^{GP}(1,Na)$ converge to Thomas-Fermi minimizer $(8\pi)^{-1/2}[\tilde{\mu}-V(x)]_+^{1/2}$ with $\tilde{\mu}$ is a suitable Lagrange multiplier. Therefore, \textit{it is natural to ask whether such a conclusion holds for the cubic-quintic nonlinear Schr\"odinger energy}, and our answer is positive, see Theorem \ref{the2} below.

\textbf{The second motivation} comes from the work  about the Ginzburg-Landau energy (see \cite{CV1993Asymptotics})
\begin{gather*}\label{GL}
	\begin{aligned}
		E^{GL}_\varepsilon(u)=\frac{1}{2}\int_{\Omega}|\nabla u|^2\,dx
		+\frac{1}{4\varepsilon^2}\int_{\Omega}\left(|u|^2-1\right)^2\,dx,   \ \ \Omega\subset \mathbb{R}^2,
	\end{aligned}
\end{gather*}
minimized in $\mathscr{H}_g=\{u\in H^1(\Omega): \int_{\Omega}|u|^2\,dx=1, u=g(x)\ \ \text{in}\ \ \partial \Omega\}$
and   the Gross-Pitaevskii energy (see \cite{Aftalion-Noris-Sourdis-CMP2015,Aftalion-ARMA-05,Ignat-RMP-06,Karali-Sourdis-ARMA2015, Aftalion-Jerrard-Royo-JFA2011,Clement-JMP2013} and references therein)
\begin{gather*}\label{GP}
	\begin{aligned}
		E^{GP}_\varepsilon(u)=\frac{1}{2}\int_{\mathbb{R}^2}|\nabla u|^2\,dx+\frac{1}{2\varepsilon^2}\int_{\mathbb{R}^2}V(x)|u|^2\,dx
		+\frac{1}{4\varepsilon^2}\int_{\mathbb{R}^2}|u|^4\,dx
	\end{aligned}
\end{gather*}
minimized in $\mathscr{H}_N=\{u\in H^1(\mathbb{R}^2): \int_{\mathbb{R}^2} V(x)|u |^2<\infty, \int_{\mathbb{R}^2}|u|^2\,dx=N\}$.  In \cite{CV1993Asymptotics}, the authors proved that proved $L^\infty$ convergence rate that $\|u_\varepsilon-u_0\|_{L^{\infty}(\Omega)}\leq{O}(\varepsilon^2)$ as $\varepsilon\to0^+$, where $u_\varepsilon$ be a minimizer for $E^{GL}_\varepsilon(\cdot)$ on $\mathcal{H}^1_g$ and $u_0$ is identified. Karali and Sourdis \cite{Karali-Sourdis-ARMA2015} studied the  ground state of a Gross-Pitaevskii energy with general trapping potential in the Thomas-Fermi limit, the ground state has irregular behavior in the form of a steep corner layer near the boundary $\partial\Omega$ of compact support of the Thomas-Fermi minimizer, they also give the $L^\infty$-convergence rates. 
 
Inspired by above, we will investigate  the $L^\infty$-convergence rate of the ground state of the cubic-quintic Schr\"{o}dinger energy, see the result of Theorem \ref{the3} below. We should mention that the methods used to estimate the  $L^\infty$-convergence rates in  the results above are  heavily reliant on the special  algebraic properties of the Gross-Pitaevskii (or Ginzburg-Landau) nonlineaity, and they do not apply to cubic-quintic Schr\"{o}dinger energy as discussed in this paper. 
 In this work, we establish a novel  energy method to obtain the $L^\infty$-convergence rate, which may represent a pioneering approach to deal with general nonlinearities.  Our method is also applicable to  the cases without potential, see \cite{lw-arxiv-b}.

\subsection{Main Results}
Based on the above motivation, we noted that the energy $e(N)$  has the scaling property 
\begin{gather}
    \begin{aligned} 
		e(N)=N \mathscr{E}(N),
    \end{aligned}
\end{gather}
where the constraint minimizing problem $\mathscr{E}(N)$ defined by 
\begin{gather}\label{amin}
    \begin{aligned} 
		\mathscr{E}(N)=\inf\left\lbrace E_N(\varphi):\varphi\in\mathscr{H}_1\right\rbrace.
    \end{aligned}
\end{gather}
The cubic-quintic Schr\"{o}dinger energy functional  $E_N(\varphi)$ is given by
\begin{gather}\label{aenergy}
    \begin{aligned} 
		{E}_N(\varphi)=\frac{1}{2}\int_{\mathbb{R}^d}\left(|\nabla \varphi|^2+V(x)|\varphi|^2\right)\,dx
		+\frac{\kappa N}{4}\int_{\mathbb{R}^d}|\varphi|^4\,dx
		+\frac{N^2}{6}\int_{\mathbb{R}^d}|\varphi|^6\,dx,
    \end{aligned}
\end{gather}
and the space $\mathscr{H}_1$ is defined as 
\begin{gather}\label{space}
    \begin{aligned} 
		\mathscr{H}_1=\left\lbrace \varphi\in H^1(\mathbb{R}^d):\ \int_{\mathbb{R}^d}V(x)|\varphi|^2\,dx<\infty, \int_{\mathbb{R}^d}|\varphi|^2\,dx=1\right\rbrace.
    \end{aligned}
\end{gather}
Indeed, it is easy to check that $\varPhi\in\mathscr{H}_N$ is a  minimizer of  $e(N)$ is equivalent with $\varphi_N(x)$ is a minimizer of the minimization problem $\mathscr{E}(N)$ and $\varphi_N(x)=\varPhi/\sqrt{N}\in\mathscr{H}_1$. In what follows, we will consider the equivalent minimization problem $\mathscr{E}(N)$. 

The crucial difference between \eqref{aenergy} and \eqref{GP} (or \eqref{GL}) is that the cubic-quintic energy functional $E_N(\cdot)$ has cubic-quintic mixed nonlinearities, while the Gross-Pitaevskii energy functional $E^{GP}_\varepsilon(\cdot)$ has only pure cubic nonlinearity. As a result, we are not able to apply Ginzburg-Landau techniques similar to the Gross-Pitaevskii functional $E^{GP}_\varepsilon(\cdot)$ to deal with this problem, which is one of the difficulties in this paper. 

For this reason, we need a new way to solve this difficulty, namely, \emph{energy methods}.  Roughly speaking, the proof of the main theorem is based on precise estimates of the energy $\mathscr{E}(N)$ and $L^p$-norm of the non-negative radially symmetric ground state $\varphi_{N}\in\mathscr{H}_1$. In fact, we shall show that
\begin{align*}\label{}
	\mathscr{E}(N)\sim N^{\frac{2p}{2d+p}}\text{ and }	\int_{\mathbb{R}^d}|\varphi_{N}|^6\,dx\sim N^{-\frac{4d}{2d+p}}\quad\text{as}\quad N\to+\infty.
\end{align*}
%where  $\varphi_{N}\in\mathscr{H}_1$ is a non-negative radially symmetric  ground state of $\mathscr{E}(N)$. 
This means that we need to consider the $L^2$-normalized scaling $w_\tau(x)=N^{\frac{d}{2d+p}}\varphi_{N}(N^{\frac{2}{2d+p}}x)$. Moreover, the energy functional $E_N(\cdot)$ can be rewritten as
\begin{gather*}\label{}
	\begin{aligned} 
		E_N(\varphi_{N})&=N^{\frac{2p}{2d+p}}\left(\frac{N^{-\frac{2p}{2d+p}}}{2}\int_{\mathbb{R}^d}V\left(N^{\frac{2}{2d+p}}x\right)|w_\tau|^2\,dx
		+\frac{1}{6}\int_{\mathbb{R}^d}|w_\tau|^6\,dx\right. \\
		&\qquad\qquad\quad\left. +\frac{N^{-\frac{2(p+2)}{2d+p}}}{2}\int_{\mathbb{R}^d}|\nabla w_\tau|^2\,dx-\frac{N^{-\frac{p}{2d+p}}}{4}\int_{\mathbb{R}^d}|w_\tau|^4\,dx\right),
	\end{aligned}
\end{gather*}
which suggests that  the last two integral above should be close to zero as $N\rightarrow+\infty$. Indeed, it is easy to know from the Sobolev embedding theorem that the cubic term  in the cubic-quintic energy functional ${E}_N(\cdot)$ is negligible as $N\to+\infty$.  A limit often considered in the literature, and called the Thomas-Fermi regime, occurs when the kinetic energy is small in front of the trapping and interaction (quintic) terms. This corresponds to the so-called Thomas-Fermi limit where the gradient and cubic terms have been dropped altogether and the functional is
\begin{gather}\label{TF1}
	\begin{aligned} 
		E^{TF}(u)=\frac{1}{2}\int_{\mathbb{R}^d}C_0|x|^p|u|^2\,dx+\frac{1}{6}\int_{\mathbb{R}^d}|u|^6\,dx,
	\end{aligned}
\end{gather}
where $C_0>0$ is a constant and given by \eqref{144}. In fact, later in the proof we will see that the limit depends only on the properties of the external potential $V(x)$ at infinity. In other words, it's independent of the shape in the bounded region of $V(x)$. The ground state energy of the Thomas-Fermi functional \eqref{TF1} with the subsidiary condition $\int_{\mathbb{R}^d}|u|^2\,dx=1$ is denoted $e^{TF}(\infty)$, that is,
\begin{align}\label{eTF1}
	e^{TF}(\infty)=\inf\left\lbrace E^{TF}(u):\  u\in L^2(\mathbb{R}^d, |x|^p\,dx)\cap L^6(\mathbb{R}^d)\ \ \text{and}\ \int_{\mathbb{R}^d}|u|^2\,dx=1\right\rbrace.
\end{align}

This also explains the seemingly strange conclusions in the following Theorems \ref{the2} and  \ref{the3}. To obtain the existence and uniqueness of the minimizer $u^{TF}$ of Thomas-Fermi energy $e^{TF}(\infty)$, one can use the arguments in \cite[Chapter 11]{anailsis}, or  see \cite{PRALieb,KRMB,Karali-Sourdis-ARMA2015} where the Thomas-Fermi limit was considered for Gross-Pitaevskii energy.  In the context of Bose-Einstein condensates, the minimizer $u^{TF}$  is known as the \emph{Thomas-Fermi approximation}.

\emph{Our first main result} concerns the limiting behavior of ground states of $\mathscr{E}(N)$ as $N\to+\infty$. Specifically, we have the following convergence theorem:
\begin{thm}\label{the2}
	Assume  that $1\leq d\leq3$, $V(x)$ satisfy \textbf{(V$_1$)}-\textbf{(V$_4$)} and $\varphi_{N}(x)$ is a non-negative radially symmetric ground state of $\mathscr{E}(N)$. Then there exists a subsequence  $\{N_k\}$  with $N_k\to+\infty$ as $k\to+\infty$, such that
	\begin{align}\label{h1}
		N_k^{\frac{d}{2d+p}}\varphi_{N_k}(N_k^{\frac{2}{2d+p}}x) \to u^{TF}(x) \ \ \text{strongly in}\  L^q(\mathbb{R}^d)\ \text{for}\ 2\leq q\leq 6,
	\end{align}	
	where $u^{TF}(x)$ is the unique minimizer of $e^{TF}(\infty)$ given by 
	\begin{gather}\label{utf1}
		\begin{aligned} 
			u^{TF}(x)=\left[\mu^{TF}-C_0|x|^p\right]^{\frac{1}{4}}_{+} \ \ \ \text{with}\ \ \ 	\mu^{TF}=\left(\frac{p}{\omega_dC_0^{-\frac{d}{p}}\mathcal{B}\left(\frac{d}{p},\frac{3}{2}\right)}\right)^{\frac{2p}{2d+p}}.
		\end{aligned}
	\end{gather}
	Here $\left[f\right]_+:=\max\{0,f\}$,  $\omega_d$ denotes the surface of the unit sphere in $\mathbb{R}^d$ and $\mathcal{B}(\cdot,\cdot)$ denotes Beta function.
	Moreover,  we have 
	\begin{align}\label{h2}
		\lim_{N\to+\infty}\frac{\mathscr{E}(N)}{N^{\frac{2p}{2d+p}}}=\frac{\omega_d(2d+p)C_0^{-\frac{d}{p}}\left(\mu^{TF}\right)^{\frac{2d+3p}{2p}}}{4pd}\mathcal{B}\left( \frac{d+p}{p},\frac{3}{2}\right),
	\end{align}
	and the following sharp vanishing rate, 
	\begin{align}\label{h3}
		\|\varphi_{N}\|_{L^{\infty}(\mathbb{R}^d)}\sim N^{-\frac{d}{2d+p}} \ \ (N\rightarrow +\infty).
	\end{align}
\end{thm} 

\begin{re}
The limit \eqref{h3} implies the vanishing phenomenon of the ground state $\varphi_{N}(x)$ at $N\to+\infty$. Similarly, the vanishing phenomenon occurs in the two-dimensional cubic-quintic nonlinear Schr\"odinger equation was studied by the authors of this paper in \cite{lw-24} as the $L^2$-mass constraint $N$ tends to the critical value $N_*$ from above.
This is very different from the purely cubic equation (with an external potential), one can see  \cite{Guo_Seiringer2014,Qwang} and the references therein. In particular, Guo and Seiringer \cite{Guo_Seiringer2014} studied the mass concentration (blow-up) properties of ground states for the purely cubic equation with an external potential as $N\nearrow N_*$ ($N<N_*$ and $N$ tends to $N_*$). The second author Wang and Feng \cite{Qwang} studied the blow-up properties of ground state solutions of the two-dimensional cubic-quintic nonlinear Schr\"odinger equation with harmonic potential, and it can also be seen that ref \cite{lw-24cv} have some similar results for other models.

\end{re}

\emph{Another main result} concerns the $L^\infty$-convergence rate of ground states $\varphi_{N}$ of the energy $\mathscr{E}(N)$ as $N\rightarrow +\infty$. 

%Our main result is the following:

\begin{thm}\label{the3}
Assume  that $1\leq d\leq3$, $V(x)$ satisfy \textbf{(V$_1$)}-\textbf{(V$_4$)} and $\varphi_{N}(x)$ is a non-negative radially symmetric ground state of $\mathscr{E}(N)$.  Let $\sigma=\min\{p-\alpha,1\}$, $\tau=N^{-\frac{2}{2d+p}}$ and $ w_\tau(x)=\tau^{-d/2}\varphi_{N}({x}/{\tau})$. Then, we have, as $\tau\to0^+(N\to+\infty)$,
    \begin{align}\label{114}
    \left| \mu_\tau-\mu^{TF}\right| \lesssim\tau^\sigma&,\\
    \left\| \Delta w_\tau(x)\right\|_{L^\infty\left( B\left( \sqrt[p]{\mu^{TF}/C_0}\right) \right) }&\lesssim\tau^{-\frac{p+10}{4}}.\label{115}
    \end{align}	
Moreover, for any $0<\epsilon\leq\sigma/2$, we have\\
(i) (The compact subset $K$ of $B(\sqrt[p]{\mu^{TF}/C_0})$) $\left\|w_\tau-u^{TF}\right\|_{L^\infty(K)}\lesssim\tau^{\sigma}$ as $\tau\to0^+$;\\
(ii) (Inside the corner layer)  $\left\|w_\tau-u^{TF}\right\|_{L^\infty\left( B\left( \sqrt[p]{\mu^{TF}/C_0}-\left( \tau|\ln \tau|\right)^\epsilon\right)\right) }\lesssim\tau^{\sigma-\epsilon}|\ln \tau|^{-\epsilon}$ as $\tau\to0^+$;\\
(iii) (Outside the corner layer) 
\begin{align}\label{117}
		w_\tau(x)\lesssim|x|^{-\frac{d-1}{2}}e^{-\frac{\beta\tau^{-\frac{p+4}{4}}}{2}|x|}\quad\text{for any }x\in \mathbb{R}^d\setminus B\left( \sqrt[p]{\mu^{TF}/C_0}+ \tau^{{p}/{2}-\epsilon}|\ln \tau|^{-\epsilon}\right).
\end{align}
Here the Lagrange multiplier $\mu_\tau$ and $\mu^{TF}$ are associated to $w_\tau(x)$ and $u^{TF}(x)$ (see  \eqref{hELw}  and \eqref{ELTF} below).
\end{thm}

\textbf{The comments to Theorem \ref{the3}}:

1. The main highlight of our paper is that we establish an novel  energy method to prove the $L^\infty$-convergence rate given in (i) and (ii) above. Our  arguments are  depending on the algebraic  properties of the nonlinearity, may be applicable to other general nonlinearities .

2. The main difference from Ginzburg-Landau  functional \cite{CV1993Asymptotics} is that the Thomas-Fermi approximation $u^{TF}$ vanishes on the boundary of the domain $B( \sqrt[p]{\mu^{TF}/C_0})$  and the square of $|\nabla u^{TF}|$ is not integrable. Thus, we will have to restrict our analysis to region $\mathcal{D}_\delta=\{x\in B( \sqrt[p]{\mu^{TF}/C_0}):\text{dist}(x,\partial B( \sqrt[p]{\mu^{TF}/C_0}))>\delta \}$ for some sufficiently small $\delta>0$.	Our method certainly allows to relax these assumptions in various ways. In fact, conditions  \textbf{(V$_2$)} and  \textbf{(V$_3$)} on the external potential $V$ can be completely unnecessary for the purposes of Theorem \ref{the2}, but they are quite important for the following result.

3. The proof of (ii) relies on the expected size of the boundary layer, namely $\left( \tau|\ln \tau|\right)^\epsilon$ for any $0<\epsilon\leq\frac{1}{2}\min\{p-\alpha,1\}$. In fact, the chosen $\left( \tau|\ln \tau|\right)^\epsilon$ may not be optimal, and the boundary layer $\delta_\tau$ only needs to satisfy $\delta_\tau\to0$ and $\delta_\tau^2/\tau^{\sigma}\to+\infty$ as $\tau\to0^+$.

\vspace{0.2cm}

\noindent
\textbf{Notation:}\\
- $\mathcal{B}(P,Q)$ denotes the Beta function with $P,Q\in\mathbb{R}^+$, i.e.,
\begin{align*}
	\mathcal{B}(P,Q)=\int_{0}^{1}x^{P-1}(1-x)^{Q-1}\,dx.
\end{align*}
- $\|\cdot\|_{s}$ denotes the $L^s(\mathbb{R}^d)$ norm for $s\geq 1$.\\
- $\tau\to0^+$ denotes $\tau>0$ such that $\tau$ tends to $0$.\\
- $\left[f\right]_+:=\max\{0,f\}$.\\  
- $\omega_N$ denotes the surface of the unit sphere in $\mathbb{R}^N$.\\
- $f^*$ denotes the symmetric-decreasing rearrangement of the function $f$.\\
- The value of positive constant $C$ is allowed to change from line to line and also in the same formula.\\
- $X\lesssim Y\ (X\gtrsim Y)$ denotes $X\leq CY\ (X\geq CY)$ for some appropriate positive constants $C$.\\
- $X\sim Y$ denotes $X\lesssim Y$ and  $Y\lesssim X$.\\

\section{Energy Estimates}\label{sec2}
Suppose that $\varphi_{N}\in\mathscr{H}_1$ is a non-negative ground state of $\mathscr{E}(N)$, then $\varphi_{N}$ satisfies the Euler-Lagrange equation
\begin{align}\label{EL}
	-\Delta  \varphi_{N}+V(x)\varphi_{N}+\kappa N \varphi_{N}^3+N^2 \varphi_{N}^{5}
	=\mu_{N}\varphi_{N},\ \ x\in\mathbb{R}^d,
\end{align}
where the corresponding Lagrange multiplier  $\mu_{N}\in\mathbb{R}$ can be computed as (multiply \eqref{EL} by $\varphi_{N}$ and integrate over $x$)
\begin{align}\label{mulam}
	\mu_N=2\mathscr{E}(N)+\frac{\kappa N}{2}\int_{\mathbb{R}^d}|\varphi_{N}|^4\,dx+\frac{2N^2}{3}\int_{\mathbb{R}^d}|\varphi_{N}|^6\,dx.
\end{align}
Moreover, we have the Pohozaev-type identity $\partial_\eta {E_N}(\eta^{\frac{d}{2}} \varphi_{N}(\eta x))|_{\eta=1}=0$ (see \cite{CIMStc} or \cite{TSPjz}), i.e.,
\begin{align}\label{PI}
	\int_{\mathbb{R}^d}|\nabla \varphi_{ N}|^2\,dx-\frac{1}{2}\int_{\mathbb{R}^d}(\nabla V\cdot x)|\varphi_{N}|^2\,dx +\frac{d\kappa N}{4}\int_{\mathbb{R}^d}|\varphi_{N}|^4\,dx+\frac{dN^2}{3}\int_{\mathbb{R}^d}|\varphi_{N}|^6\,dx=0.
\end{align}

We recall from \cite{CMPwein} the well-known Gagliardo-Nirenberg inequality with the best constant:
\begin{gather}\label{GN}
	\begin{aligned} 
		\int_{\mathbb{R}^d}|u|^4\,dx\leq\frac{2}{\|Q\|_2^{2}}\left( \int_{\mathbb{R}^d}|\nabla  u|^2\,dx\right)^{\frac{d}{2}}\left(\int_{\mathbb{R}^d}|u|^2\,dx\right)^{\frac{4-d}{2}},\quad u\in H^1(\mathbb{R}^d),
	\end{aligned}
\end{gather}
which can be equality only when $u=Q(|x|)$, where  $Q$ is the unique  positive radial ground state solution of
the following scalar-field equation
\begin{gather}\label{Qequ}
	\begin{aligned} 
		-\frac{d}{2}\Delta u+\frac{4-d}{2}u-u^3=0,\quad u\in H^1(\mathbb{R}^d).
	\end{aligned}
\end{gather}
%Recall from \cite{CMPwein,kmk1989uniqueness,wm1996Progress} that, up to translations, \eqref{Qequ} possesses a unique positive solution $Q > 0$,
%which must be radially symmetric. 
Also, one can derive that a unique solution $Q$ of \eqref{Qequ} satisfies
\begin{gather}\label{Qgj}
	\begin{aligned}
		\int_{\mathbb{R}^d}|\nabla Q|^2\,dx
		=\int_{\mathbb{R}^d}|Q|^2\,dx
		=\frac{1}{2}\int_{\mathbb{R}^d}|Q|^4\,dx.
	\end{aligned}
\end{gather}
Note also  that $Q(|x|)$ has the following exponential 
decay in the following sense:
\begin{gather}\label{decay}
	\begin{aligned} 
		Q(|x|), |\nabla Q(|x|)|=O(|x|^{-\frac{1}{2}}e^{-\nu|x|})\quad\text{as}\ |x|\to\infty,
	\end{aligned}
\end{gather}
where $\nu$ is a positive constant. 

Next, we prove the following energy estimates. 
\begin{lem}\label{lem2}
	Assume that $1\leq d\leq3$, $V(x)$ satisfies \textbf{(V$_1$)} and \textbf{(V$_4$)}, we then  have the following energy estimates
	\begin{align}\label{37}
		\mathscr{E}(N)\sim N^{\frac{2p}{2d+p}}\quad\text{as}\quad N\to+\infty.
	\end{align}
\end{lem}
\begin{pf} We start with the upper bound estimate on the energy $\mathscr{E}(N)$ as $ N\to+\infty$. We deduce from \textbf{(V$_4$)} that there exists a constant $R>1$ such that
\begin{gather}\label{32}
    \begin{aligned} 
		\frac{1}{R}|x|^p\leq V(x)\leq R|x|^p, \quad \text{if }|x|\geq R.
    \end{aligned}
\end{gather}
 
Let $Q(|x|)$ is the unique radially symmetric positive solution of \eqref{Qequ}, for any constant $\eta>0$, choose a trial function
\begin{gather*}
    \begin{aligned} 
		Q_{\eta}(x)=\frac{\eta^{\frac{d}{2}} Q(\eta x)}{\|Q\|_{2}},
    \end{aligned}
\end{gather*}
thus, we fined that $Q_{\eta}\in\mathscr{H}_1$.	Moreover, using by \eqref{Qgj} and \eqref{decay}, let $\eta>0$ is sufficiently small, then 
\begin{gather*}
    \begin{aligned} 
		{E}_N(Q_{\eta})
		&\leq\frac{1}{2}\int_{\mathbb{R}^d}|\nabla Q_{\eta}|^2\,dx
		+\frac{1}{2}\int_{B(0,R)}V(x)|Q_{\eta}|^2\,dx
		+\frac{1}{2}\int_{\mathbb{R}^d\setminus B(0,R)}V(x)|Q_{\eta}|^2\,dx\\
		&\quad+\frac{N}{4}\int_{\mathbb{R}^d}|Q_{\eta}|^4\,dx
		+\frac{N^2}{6}\int_{\mathbb{R}^d}|Q_{\eta}|^6\,dx\\
		&\leq\frac{1}{2}\int_{\mathbb{R}^d}|\nabla Q_{\eta}|^2\,dx
		+\frac{1}{2}\int_{B(0,R)}M|Q_{\eta}|^2\,dx
		+\frac{R}{2}\int_{\mathbb{R}^d\setminus B(0,R)}|x|^p|Q_{\eta}|^2\,dx\\
		&\quad+\frac{N}{4}\int_{\mathbb{R}^d}|Q_{\eta}|^4\,dx
		+\frac{N^2}{6}\int_{\mathbb{R}^d}|Q_{\eta}|^6\,dx\\
		&\leq\frac{M}{2}+\frac{\eta^2}{2}
		+\frac{R\eta^{-p}}{2\|Q\|_{L^2}^2}\int_{\mathbb{R}^d}|x|^p|Q|^2\,dx	
		+\frac{N\eta^{d}}{2\|Q\|_{L^2}^2}
		+\frac{N^2\eta^{2d}}{6\|Q\|_{L^2}^6}\int_{\mathbb{R}^d}|Q|^6\,dx\\
		&=:\frac{M}{2}+C_1\eta^2+C_2\eta^{-p}+C_3N\eta^{d}+C_4N^2\eta^{2d},
    \end{aligned}
\end{gather*}	 
where $M:=\max\{V(x):x\in B(0,R)\}$ and   constants $C_i>0\ (i=1,2,3,4)$, independent of $\eta$ and $N$.  Taking  $\eta=N^{-\frac{2}{2d+p}}$, we deduce that 
\begin{gather*}
    \begin{aligned} 
		{E}_N(Q_{\eta})\leq\frac{M}{2}+ C_1N^{-\frac{4}{2d+p}}+C_2N^{\frac{2p}{2d+p}}+C_4N^{\frac{p}{2d+p}}+C_4N^{\frac{2p}{2d+p}}.
    \end{aligned}
\end{gather*}	
Since $V\in C^1$, then $M$ exists and is bounded in  $B(0,R)$.	Therefore, we have 
\begin{align}\label{219}
    \mathscr{E}(N)\lesssim N^{\frac{2p}{2d+p}}\qquad\text{as}\ N\to+\infty.
\end{align} 
	
Next, we shall establish the lower bound estimate of $\mathscr{E}(N)$ as $N\to+\infty$.  Let a constant $\gamma>R^{p-1}$.	On the one hand, for any $v\in\mathscr{H}_1$, by direct calculation, we obtain 
\begin{gather}\label{5.8}
    \begin{aligned} 
		&\quad\frac{1}{2}\int_{\mathbb{R}^d}V(x)|v(x)|^2\,dx
		+\frac{N^2}{6}\int_{\mathbb{R}^d}|v(x)|^6\,dx\\
		&\geq\frac{1}{2}\int_{\mathbb{R}^d\setminus B(0,R)}V(x)|v(x)|^2\,dx
		+\frac{N^2}{6}\int_{\mathbb{R}^d}|v(x)|^6\,dx\\
		&\geq\frac{\gamma}{2}+\frac{1}{2}\int_{\mathbb{R}^d\setminus B(0,R)}\left( V(x)-\gamma\right) |v(x)|^2\,dx-\frac{\gamma}{2}\int_{B(0,R)} |v(x)|^2\,dx
		+\frac{N^2}{6}\int_{\mathbb{R}^d}|v(x)|^6\,dx\\
		&\geq\frac{\gamma}{2}-\frac{1}{2}\int_{\mathbb{R}^d\setminus B(0,R)}\left[\gamma-V(x)\right]_+|v(x)|^2\,dx
		+\frac{N^2}{6}\int_{\mathbb{R}^d\setminus B(0,R)}|v(x)|^6\,dx\\
		&\quad-\frac{\gamma}{2}\int_{B(0,R)} |v(x)|^2\,dx
		+\frac{N^2}{6}\int_{B(0,R)}|v(x)|^6\,dx\\
		&\geq\frac{\gamma}{2}-\frac{1}{2}\int_{\mathbb{R}^d\setminus B(0,R)}\left[\gamma-R^{-1}|x|^p\right]_+|v(x)|^2\,dx
		+\frac{N^2}{6}\int_{\mathbb{R}^d\setminus B(0,R)}|v(x)|^6\,dx\\
		&\quad-\frac{\gamma}{2}\int_{B(0,R)} |v(x)|^2\,dx
		+\frac{N^2}{6}\int_{B(0,R)}|v(x)|^6\,dx,
    \end{aligned}
\end{gather}
where  $\left[f\right]_+:=\max\{0,f\}$. Using by Young's inequality, we have 
\begin{align} \label{5.91}
	\frac{1}{2}\left[\gamma-R^{-1}|x|^p\right]_+|v(x)|^2\leq\frac{N^2}{6}|v(x)|^6+\frac{1}{3N}\left[\gamma-R^{-1}|x|^p\right]_+^{\frac{3}{2}},
\end{align}
and
\begin{align} \label{5.92}
	\frac{\gamma}{2}|v(x)|^2\leq\frac{N^2}{6}|v(x)|^6+\frac{1}{3N}\gamma^{\frac{3}{2}}.
\end{align}
Note that
\begin{gather}\label{5.10123}
    \begin{aligned} 
		\int_{\mathbb{R}^d\setminus B(0,R)}\left[\gamma-R^{-1}|x|^p\right]_+^{\frac{3}{2}}\,dx
		&=\int_{R\leq |x|\leq{\sqrt[p]{R\gamma}}}(\gamma-R^{-1}|x|^p)^{\frac{3}{2}}\,dx\\
		&=\int_R^{\sqrt[p]{R\gamma}}\,dr\int_{\partial B(0,r)}(\gamma-R^{-1}r^p)^{\frac{3}{2}}\,dS\\
		%=&\omega_d\int_R^{\sqrt[p]{R\gamma}}(\gamma-R^{-1}r^p)^{\frac{3}{2}}r^{d-1}\,dr\\
		&\leq\omega_d\int_0^{\sqrt[p]{R\gamma}}(\gamma-R^{-1}r^p)^{\frac{3}{2}}r^{d-1}\,dr\\
		&\leq\omega_d\int_0^{\sqrt[p]{R\gamma}}\gamma^{\frac{3}{2}}r^{d-1}\,dr\\
		&\leq \frac{\omega_dR^{\frac{d}{p}}}{d}\gamma^{\frac{2d+3p}{2p}},
    \end{aligned}
\end{gather}
where $\omega_d$ denotes the surface of the unit sphere in $\mathbb{R}^d$. Thus, we deduce from \eqref{5.8}--\eqref{5.10123}  that
	\begin{gather*}%\label{5.10}
		\begin{aligned} 
			&\quad\frac{1}{2}\int_{\mathbb{R}^d}V(x)|v(x)|^2\,dx
			+\frac{N^2}{6}\int_{\mathbb{R}^d}|v(x)|^6\,dx\\
			&\geq\frac{\gamma}{2}-\frac{\omega_dR^d}{3dN}\gamma^{\frac{3}{2}}
			-\frac{\omega_dR^{\frac{d}{p}}}{3dN}\gamma^{\frac{2d+3p}{2p}}\\
			&=\gamma\left(\frac{1}{2}-\frac{\omega_dR^d}{3dN}\gamma^{\frac{1}{2}}
			-\frac{\omega_dR^{\frac{d}{p}}}{3dN}\gamma^{\frac{2d+p}{2p}} \right) .
		\end{aligned}
	\end{gather*}
	Taking $\gamma=CN^{\frac{2p}{2d+p}}$ ($ >R^{p-1}$ for large enough $N>0$), 
	$C>0$ is a small constant such that $0<\omega_dR^{\frac{d}{p}}C^{\frac{2d+p}{2p}}/(3d)<\frac{1}{4}$, we then have
	\begin{gather}\label{214}
		\begin{aligned} 
			&\quad\frac{1}{2}\int_{\mathbb{R}^d}V(x)|v(x)|^2\,dx
			+\frac{N^2}{6}\int_{\mathbb{R}^d}|v(x)|^6\,dx\\
			&\gtrsim N^{\frac{2p}{2d+p}}\left(\frac{1}{2}-\frac{\omega_dR^d}{3d}N^{-\frac{2d}{2d+p}}
			-\frac{\omega_dR^{\frac{d}{p}}C^{\frac{2d+p}{2p}}}{3d} \right)\\
			&\gtrsim N^{\frac{2p}{2d+p}}\quad \text{ as }N\to+\infty.
		\end{aligned}
	\end{gather}

On the other hand, let $\varphi_N(x)$ be a non-negative ground state of $\mathscr{E}(N)$ with $\|\varphi_N\|^2_{2}=1$, we deduce from  \eqref{214} that
	\begin{gather}\label{}
		\begin{aligned} 
			\mathscr{E}(N)%&=\frac{1}{2}\int_{\mathbb{R}^d}|\nabla\varphi_N|^2\,dx
			%+\frac{1}{2}\int_{\mathbb{R}^d}V(x)|\varphi_N|^2\,dx
			%+\frac{\kappa N}{4}\int_{\mathbb{R}^d}|\varphi_N|^4\,dx
			%+\frac{N^2}{6}\int_{\mathbb{R}^d}|\varphi_N|^6\,dx\\
			&\geq\frac{1}{2}\int_{\mathbb{R}^d}V(x)|\varphi_N|^2\,dx
			+\frac{N^2}{6}\int_{\mathbb{R}^d}|\varphi_N|^6\,dx-\frac{N}{4}\int_{\mathbb{R}^d}|\varphi_N|^4\,dx\\
			&\gtrsim N^{\frac{2p}{2d+p}}-\frac{N}{4}\int_{\mathbb{R}^d}|\varphi_N|^4\,dx.
		\end{aligned}
	\end{gather}
 
Now,  we \textbf{claim that} 
\begin{gather}\label{0}
   \begin{aligned} 
   \int_{\mathbb{R}^d}|\varphi_N|^4\,dx\lesssim N^{-\frac{2d}{2d+p}}\quad \text{ as }N\to+\infty.
   \end{aligned}
\end{gather}	
Indeed, using by H\"oder inequality, we get
	\begin{gather}\label{217}
		\begin{aligned} 
			\int_{\mathbb{R}^d}|\varphi_N|^4\,dx
			\leq	\left(\int_{\mathbb{R}^d}|\varphi_{N}|^2\,dx\right)^{\frac{1}{2}}\left(\int_{\mathbb{R}^d}|\varphi_{N}|^6\,dx\right)^{\frac{1}{2}}
			=\left(\int_{\mathbb{R}^d}|\varphi_{N}|^6\,dx\right)^{\frac{1}{2}},
		\end{aligned}
	\end{gather}
where using the $\varphi_{N}\in\mathscr{H}_1$. Note that, by \eqref{219}, we have
	\begin{gather*}
		\begin{aligned} 
			\frac{N^2}{6}\int_{\mathbb{R}^d}|\varphi_{N}|^6\,dx
			-\frac{N}{4}\int_{\mathbb{R}^d}|\varphi_N|^4\,dx\leq \mathscr{E}(N)\lesssim N^{\frac{2p}{2d+p}}.
		\end{aligned}
	\end{gather*}	
Therefore, it follow that
	\begin{gather}\label{218}
		\begin{aligned} 
			\int_{\mathbb{R}^d}|\varphi_{N}|^6\,dx\lesssim N^{-\frac{4d}{2d+p}}+{N^{-1}}\int_{\mathbb{R}^d}|\varphi_N|^4\,dx.
		\end{aligned}
	\end{gather}	
	Combining \eqref{217}, we obtain
	\begin{gather*}
		\begin{aligned} 
			\int_{\mathbb{R}^d}|\varphi_N|^4\,dx
			\lesssim\left(N^{-\frac{4d}{2d+p}}+{N^{-1}}\int_{\mathbb{R}^d}|\varphi_N|^4\,dx\right)^{\frac{1}{2}}.
		\end{aligned}
	\end{gather*}
	Multiply both sides of the above inequality by   $N^{\frac{2d}{2d+p}}$, we get
	\begin{gather*}
		\begin{aligned} 
			N^{\frac{2d}{2d+p}}\int_{\mathbb{R}^d}|\varphi_N|^4\,dx
			\lesssim\left(1+N^{\frac{2d-p}{2d+p}}\int_{\mathbb{R}^d}|\varphi_N|^4\,dx\right)^{\frac{1}{2}}
			\leq1+N^{\frac{2d-p}{2d+p}}\int_{\mathbb{R}^d}|\varphi_N|^4\,dx,
		\end{aligned}
	\end{gather*}
which implies that for $N$ large enough,
	\begin{gather}\label{zj11}
		\begin{aligned}
  N^{\frac{2d}{2d+p}}\int_{\mathbb{R}^d}|\varphi_N|^4\,dx\lesssim 1.
		\end{aligned}
	\end{gather}	
Furthermore we see from \eqref{zj11} that \eqref{0} holds. Thus, we obtain 
	\begin{gather*}
		\begin{aligned} 
			\mathscr{E}(N)\gtrsim N^{\frac{2p}{2d+p}}\quad \text{ as }N\to+\infty.
		\end{aligned}
	\end{gather*}
This estimate and \eqref{219} imply that this completes the proof of this lemma.
\end{pf}

\begin{lem}\label{lemi6}
	Assume that $1\leq d\leq3$, $V(x)$ satisfies \textbf{(V$_1$)} and \textbf{(V$_4$)}, and $\varphi_{N}(x)$ is a non-negative  ground state of $\mathscr{E}(N)$, we then have
	\begin{align}\label{lem6}
	\int_{\mathbb{R}^d}|\varphi_{N}|^6\,dx\sim N^{-\frac{4d}{2d+p}}\quad\text{as }N\to+\infty.
	\end{align}
\end{lem}
\begin{pf} It follows from \eqref{218} and \eqref{zj11} that
\begin{gather*}
    \begin{aligned} 
    \int_{\mathbb{R}^d}|\varphi_{N}|^6\,dx&\lesssim N^{-\frac{4d}{2d+p}}-\kappa{N^{-1}}\int_{\mathbb{R}^d}|\varphi_{N}|^4\,dx\lesssim N^{-\frac{4d}{2d+p}}+N^{-\frac{10d+p}{2d+p}}.
    \end{aligned}
\end{gather*}	
Thus, the upper bound in \eqref{lem6} follows immediately from $N\to+\infty$.
	
Next, we prove the lower bound in \eqref{lem6}. Let $\tilde{N}>0$ and $\theta>1$ such that $\tilde{N}=\theta N$, we have 
	\begin{gather*}
		\begin{aligned} 
			\mathscr{E}(\tilde{N})\leq E_{\tilde{N}}(\varphi_{N})&=\mathscr{E}(N)+\frac{\kappa(\tilde{N}-N)}{4}\int_{\mathbb{R}^d}|\varphi_{N}|^4\,dx+\frac{\tilde{N}^2-N^2}{6}\int_{\mathbb{R}^d}|\varphi_{N}|^6\,dx\\
			&\leq \mathscr{E}(N)+\frac{\tilde{N}-N}{4}\left( \int_{\mathbb{R}^d}|\varphi_{N}|^6\,dx\right)^{\frac{1}{2}}+\frac{\tilde{N}^2-N^2}{6}\int_{\mathbb{R}^d}|\varphi_{N}|^6\,dx\\
			&\leq \mathscr{E}(N)+\frac{C(\theta-1)}{4}N^{\frac{p}{2d+p}}+\frac{\tilde{N}^2-N^2}{6}\int_{\mathbb{R}^d}|\varphi_{N}|^6\,dx.
		\end{aligned}
	\end{gather*}
	Hence, by Lemma \ref{lem2}, we obtain
	\begin{gather*}
		\begin{aligned} 
			\int_{\mathbb{R}^d}|\varphi_{N}|^6\,dx
			&\geq\frac{C_1\tilde{N}^{\frac{2p}{2d+p}}-C_2N^{\frac{2p}{2d+p}}}{\tilde{N}^2-N^2}-\frac{C(\theta-1)}{4\tilde{N}^2-N^2}N^{\frac{p}{2d+p}}\\
			&=\frac{N^{-{\frac{4d}{2d+p}}}(C_1\theta^{\frac{2p}{2d+p}}-C_2)}{\theta^2-1}-\frac{C}{4(\theta+1)}N^{-\frac{4d+p}{2d+p}}.
		\end{aligned}
	\end{gather*}
	Taking $\theta$ large enough such that $C_1\theta^{\frac{2p}{2d+p}}-C_2>0$, we can get 
	$$\int_{\mathbb{R}^d}|\varphi_{N}|^6\,dx\gtrsim N^{-\frac{4d}{2d+p}}\quad\text{as }N\to+\infty.$$
	This completes the proof of
	the lemma.
\end{pf}

\begin{re}\label{re31}
	Assume that $1\leq d\leq3$, $V(x)$ satisfies \textbf{(V$_1$)} and \textbf{(V$_4$)}, and $\varphi_{N}(x)$ is a non-negative  ground state of $\mathscr{E}(N)$, we then have
	\begin{align}\label{22}
		\int_{\mathbb{R}^d}V(x)|\varphi_{N}|^2\,dx\lesssim N^{\frac{2p}{2d+p}}\quad\text{as}\quad N\to+\infty.
	\end{align}
  Indeed, by H\"older inequality and Young inequality, we have
\begin{gather*}
	\begin{aligned} 
		\int_{\mathbb{R}^d}|\varphi_{N}|^4\,dx\leq\left( \int_{\mathbb{R}^d}|\varphi_{N}|^6\,dx\right)^{\frac{1}{2}}\leq\frac{2N}{3}\int_{\mathbb{R}^d}|\varphi_{N}|^6\,dx+\frac{3}{8N}.
	\end{aligned}
\end{gather*}
Hence, we have
\begin{gather*}
	\begin{aligned}
		{E_N}(\varphi_N)&\geq\frac{1}{2}\int_{\mathbb{R}^d}\left( |\nabla \varphi_N|^2+V(x)|\varphi_N|^2\right) \,dx-\frac{N}{4}\int_{\mathbb{R}^d}|\varphi_N|^4\,dx+\frac{N^2}{6}\int_{\mathbb{R}^d}|\varphi_N|^6\,dx\\
		&\geq\frac{1}{2}\int_{\mathbb{R}^d}\left( |\nabla \varphi_N|^2+V(x)|\varphi_N|^2\right) \,dx-\frac{3}{32}.
	\end{aligned}
\end{gather*}
This estimate, \eqref{37} imply finally \eqref{22} holds.
\end{re}

\section{Thomas-Fermi Approximation}

Suppose $\varphi_{N}\in\mathscr{H}_1$ be a non-negative radially symmetric ground state of $\mathscr{E}(N)$. For simplicity, we denote $\tau=N^{-\frac{2}{2d+p}}>0$ in the following proof. Define the following $L^2$-normalized function 
\begin{align}\label{hscaling}
w_\tau(x):= \tau^{-\frac{d}{2}}\varphi_{N}\left( \frac{x}{\tau}\right) .
\end{align}
Note that $\tau\to0^+$ as $N\to+\infty$. Moreover, we then have
\begin{gather}\label{h0310}
\begin{aligned} 
	E_N(\varphi_{N})=\mathscr{E}(N)=:\frac{1}{\tau^{p}}e(\tau),
\end{aligned}
\end{gather}	
where the  minimizing problem $e(\tau)$ be given by
\begin{gather}\label{hEtau0}
\begin{aligned} 
	e(\tau):=\inf\left\lbrace I_{\tau}(u):\  u\in\mathscr{H}_1\right\rbrace
\end{aligned}
\end{gather}
and the energy function $I_{\tau}(\cdot)$ be given by
\begin{gather*}
\begin{aligned}
	I_{\tau}(u):=\frac{\tau^{p+2}}{2}\int_{\mathbb{R}^d}|\nabla u|^2\,dx +\frac{\tau^{p}}{2}\int_{\mathbb{R}^d}V(\tau^{-1}x)|u|^2\,dx +\frac{\kappa\tau^{\frac{p}{2}}}{4}\int_{\mathbb{R}^d}|u|^4\,dx+\frac{1}{6}\int_{\mathbb{R}^d}|u|^6\,dx.
\end{aligned}
\end{gather*}
Thus, $w_\tau(x)$ is a minimizer of  \eqref{hEtau0}, i.e.,
\begin{gather*}
\begin{aligned} 
	I_{\tau}(w_\tau)=e(\tau).
\end{aligned}
\end{gather*}	
Moreover, we have the Pohozaev-type identity $\partial_\eta I_{\tau}(\eta^{d/2} w_\tau(\eta x))|_{\eta=1}=0$, i.e.,
\begin{align}\label{hIPI}
{\tau^{p+2}}\int_{\mathbb{R}^d}|\nabla w_\tau|^2\,dx +\frac{\kappa\tau^{\frac{p}{2}} d}{4}\int_{\mathbb{R}^d}|w_\tau|^4\,dx
=\frac{\tau^{p-1}}{2}\int_{\mathbb{R}^d}(\nabla V(\tau^{-1}x)\cdot x)|w_\tau|^2\,dx-\frac{d}{3}\int_{\mathbb{R}^d}|w_\tau|^6\,dx,
\end{align}
and $w_\tau(x)$ satisfies the following Euler-Lagrange equation 
\begin{gather}\label{hELw}
\begin{aligned} 
	-\tau^{p+2}\Delta  w_\tau(x)+\tau^{p}V(\tau^{-1}x)w_\tau(x)+\kappa\tau^{\frac{p}{2}} w_\tau^3(x)+w_\tau^5(x)
	=\mu_\tau w_\tau(x)\ \ \text{in}\ \mathbb{R}^d,
\end{aligned}
\end{gather}
where $\mu_\tau\in\mathbb{R}$ is a suitable Lagrange multiplier associated to $w_\tau(x)$.
The corresponding Lagrange multiplier $\mu_\tau$ can be computed
as (multiply \eqref{hELw} by $w_\tau$ and integrate over $x$)
\begin{gather}\label{hwmu}
\begin{aligned} 
	\mu_\tau=2e(\tau)+\frac{\kappa\tau^{\frac{p}{2}}}{2}\int_{\mathbb{R}^d}|w_\tau|^4\,dx+\frac{2}{3}\int_{\mathbb{R}^d}|w_\tau|^6\,dx,
\end{aligned}
\end{gather}
Form \eqref{zj11}, \eqref{lem6} and \eqref{22}, then for $\tau=N^{-\frac{2}{2d+p}}$ small enough we have
\begin{align}\label{410}
\int_{\mathbb{R}^d}|w_\tau|^6\,dx\sim1, \quad\int_{\mathbb{R}^d}\tau^{p}V(\tau^{-1}x)|w_\tau|^2\,dx\lesssim1\quad\text{and}\quad\int_{\mathbb{R}^d}|w_\tau|^4\,dx\lesssim1.
\end{align}

Suppose that $u^{TF}$ is the unique non-negative minimizer of \eqref{eTF1} with $\|u^{TF}\|_{2}^2=1$. Thus $u^{TF}(x)$ satisfies the following Euler-Lagrange equation 
\begin{gather}\label{ELTF}
\begin{aligned} 
	\mu^{TF}u^{TF}(x)=C_0|x|^pu^{TF}(x)+|u^{TF}(x)|^{4}u^{TF}(x)\ \ \text{in }\mathbb{R}^d,
\end{aligned}
\end{gather}
where $\mu^{TF}\in\mathbb{R}^+$ is a suitable Lagrange multiplier associated to $u^{TF}(x)$. Moreover, we have the Pohozaev-type identity $\partial_\eta E^{TF}(\eta^{d/2} u^{TF}(\eta x))|_{\eta=1}=0$, i.e.,
\begin{align}\label{hTFPI}
\frac{pC_0}{2}\int_{\mathbb{R}^d}|x|^p|u^{TF}|^2\,dx=\frac{d}{3}\int_{\mathbb{R}^d}|u^{TF}|^6\,dx.
\end{align}
The corresponding Lagrange multiplier $\mu^{TF}$ can be computed
as (multiply \eqref{ELTF} by $u^{TF}$ and integrate over $x$)
\begin{gather}\label{tfmu}
\begin{aligned} 
	\mu^{TF}
	=2e^{TF}(\infty)+\frac{2}{3}\int_{\mathbb{R}^d}|u^{TF}|^6\,dx,
\end{aligned}
\end{gather}
Applying the non-negativity of $u^{TF}(x)$, by \eqref{ELTF}, we have
\begin{gather}\label{utf}
\begin{aligned}
	u^{TF}(x)=
	\begin{cases}
		\left(\mu^{TF}-C_0|x|^p\right)^{\frac{1}{4}}, \ \  & |x|<\sqrt[p]{\mu^{TF}/C_0}, \\
		\qquad0, \ \ & \text{otherwise}.
	\end{cases}
\end{aligned}
\end{gather}
Note that
\begin{gather}\label{6133}
\begin{aligned} 
	1=\int_{\mathbb{R}^d}|u^{TF}|^2\,dx=&\int_{|x|\leq\sqrt[p]{\mu^{TF}/C_0}}(\mu^{TF}-C_0|x|^p)^{\frac{1}{2}}\,dx\\
	=&\int_{0}^{\sqrt[p]{\mu^{TF}/C_0}}\,dr\int_{\partial B(r)}(\mu^{TF}-C_0r^p)^{\frac{1}{2}}\,dS\\
	=&\omega_d\int_{0}^{\sqrt[p]{\mu^{TF}/C_0}}(\mu^{TF}-C_0r^p)^{\frac{1}{2}}r^{d-1}\,dr\\
	%=&\frac{\omega_d}{p}\int_{0}^{\mu^{TF}/C_0}(\mu^{TF}-C_0t)^{\frac{1}{2}}t^{\frac{d-p}{p}}\,dt\\
	=&\frac{\omega_dC_0^{-\frac{d}{p}}}{p}\int_{0}^{\mu^{TF}}(\mu^{TF}-t)^{\frac{1}{2}}t^{\frac{d-p}{p}}\,dt\\
	=&\frac{\omega_dC_0^{-\frac{d}{p}}\left(\mu^{TF}\right)^{\frac{2d+p}{2p}}}{p}\int_{0}^{\mu^{TF}}\left((1-\frac{t}{\mu^{TF}}\right)^{\frac{1}{2}}\left(\frac{t}{\mu^{TF}}\right)^{\frac{d-p}{p}}\,d\left(\frac{t}{\mu^{TF}} \right)\\
	=&\frac{\omega_dC_0^{-\frac{d}{p}}\left(\mu^{TF}\right)^{\frac{2d+p}{2p}}}{p}\mathcal{B}\left(\frac{d}{p},\frac{3}{2} \right),
\end{aligned}
\end{gather}
where $\omega_d$ denotes the surface of the unit sphere in $\mathbb{R}^d$ and $\mathcal{B}(\cdot,\cdot)$ denotes the Beta function. This is implies that
\begin{gather}\label{mutf}
\begin{aligned} 
	\mu^{TF}=\left(\frac{p}{\omega_dC_0^{-\frac{d}{p}}\mathcal{B}\left(\frac{d}{p},\frac{3}{2}\right)}\right)^{\frac{2p}{2d+p}}.
\end{aligned}
\end{gather}
Similar to \eqref{6133}, by direct calculation
\begin{gather}\label{N43}
\begin{aligned} 
	\int_{\mathbb{R}^d}|\nabla u^{TF}|^2\,dx
	&=\frac{p^2C_0^2}{16} \int_{|x|\leq\sqrt[p]{\mu^{TF}/C_0}}(\mu^{TF}-C_0|x|^p)^{-\frac{3}{2}}|x|^{2(p-1)}\,dx\\
	&=\frac{p^2C_0^2}{16}\int_0^{\sqrt[p]{\mu^{TF}/C_0}}\,dr\int_{\partial B(r)}(\mu^{TF}-C_0r^p)^{-\frac{3}{2}}r^{2(p-1)}\,dS\\
	&=\frac{p^2C_0^2\omega_d}{16}\int_0^{\sqrt[p]{\mu^{TF}/C_0}}(\mu^{TF}-C_0r^p)^{-\frac{3}{2}}r^{d+2p-3}\,dr\\
	%=&\frac{pC_0^2\omega_d}{16}\int_{0}^{\mu^{TF}/C_0}(\mu^{TF}-C_0t)^{-\frac{3}{2}}t^{\frac{d+p-2}{p}}\,dt\\	=&\frac{pC_0^{\frac{2-d}{p}}\omega_d}{16}\int_{0}^{\mu^{TF}}(\mu^{TF}-t)^{-\frac{3}{2}}t^{\frac{d+p-2}{p}}\,dt\\
	&=\frac{pC_0^{\frac{2-d}{p}}\omega_d(\mu^{TF})^{\frac{2d+p-4}{2p}}}{16}\int_{0}^{1}\left(1-\frac{r}{\mu^{TF}}\right)^{-\frac{3}{2}}{\left(\frac{r}{\mu^{TF}}\right)^{\frac{d+p-2}{p}}}\,d\left(\frac{r}{\mu^{TF}}\right)\\
	&=\frac{pC_0^{\frac{2-d}{p}}\omega_d(\mu^{TF})^{\frac{2d+p-4}{2p}}}{16}\mathcal{B}\left(\frac{d+2p-2}{p},-\frac{1}{2} \right).
\end{aligned}
\end{gather}
Thus,  $\nabla u^{TF}$ is not square-integrable near $|x|=\sqrt[p]{\mu^{TF}/C_0}$.

\subsection{The Gradient Estimate}
The purpose of this subsection is to discuss the blow-up estimate of the kinetic energy. This is an important part of our proof below in this paper.

\begin{lem}\label{lem52}
Assume that $1\leq d\leq3$ and $V(x)$ satisfies \textbf{(V$_1$)}--\textbf{(V$_4$)}.	Let $\sigma:=\min\{p-\alpha,1\}$, %$p\geq2$ and $0\leq\alpha<p$, 
thus, as $\tau\to0^+$, we have
\begin{gather}\label{hclaim}
	\begin{aligned} 
		\left| e(\tau)-e^{TF}(\infty)\right| \leq O(\tau^\sigma)\ \ \ \text{as }\tau\to0^+.
	\end{aligned}
\end{gather}
Moreover, if $w_\tau$ is a minimizer of $e(\tau)$, then there holds that
\begin{gather}\label{4241}
	\begin{aligned} 
		\int_{\mathbb{R}^d}|\nabla w_\tau|^2\,dx\leq O(\tau^{\sigma-p-2})\ \ \ \text{as }\tau\to0^+.
	\end{aligned}
\end{gather}
\end{lem}

\begin{pf} 
Let $0<\tau\ll1$ and $\eta(x)\in C^{\infty}_0(\mathbb{R}^d)$ be a non-negative smooth cut-off function defined by
\begin{gather}
	\begin{cases}\label{le03.1.1.17}
		\eta_\tau(x)=1
		& \ \   \text{for }\  |x|\leq\sqrt[p]{\mu^{TF}/C_0}-2\tau,\\
		0\leq\eta_\tau(x)\leq1
		& \ \   \text{for }\  \sqrt[p]{\mu^{TF}/C_0}-2\tau<|x|<\sqrt[p]{\mu^{TF}/C_0}-\tau,\\
		\eta_\tau(x)=0
		&  \ \   \text{for }\ |x|\geq\sqrt[p]{\mu^{TF}/C_0}-\tau.
	\end{cases}
\end{gather}
Without loss of generality, we may assume that $|\nabla\eta_\tau(x)|\leq{M}/{\tau}$, where $M > 0$ is a positive
constant independent of $\tau$. Consider the following  function
\begin{gather}\label{521}
	\begin{aligned} 
		u^{TF}_\tau(x):=A_\tau\eta_\tau(x)u^{TF}(x),\ \ \ x\in\mathbb{R}^d,
	\end{aligned}
\end{gather}
where $A_\tau$ is chosen such that $\int_{\mathbb{R}^d}|u^{TF}_\tau|^2\,dx=1$. Define
\begin{gather}\label{Omega}
	\begin{aligned} 
		\Omega^1_\tau&:=\left\lbrace x\in\mathbb{R}^d:\sqrt[p]{\mu^{TF}/C_0}-\tau<|x|\leq\sqrt[p]{\mu^{TF}/C_0} \right\rbrace;\\
		\Omega^2_\tau&:=\left\lbrace x\in\mathbb{R}^d:\sqrt[p]{\mu^{TF}/C_0}-2\tau<|x|\leq\sqrt[p]{\mu^{TF}/C_0}-\tau \right\rbrace.
	\end{aligned}
\end{gather}
Obviously, $u^{TF}_\tau\in C^{\infty}_0(\mathbb{R}^d)$ and $u^{TF}_\tau(x)\to u^{TF}(x)$ strongly in $L^r(\mathbb{R}^d)\cap L^2(\mathbb{R}^d,d\mu)$ with $d\mu=|x|^pdx$, where  $2\leq r\leq 6$ and $p\geq2$. Thus, we deduce that
\begin{gather}\label{hETFo1}
	\begin{aligned} 
		E^{TF}(u^{TF}_\tau)=E^{TF}(u^{TF})+o(1).
	\end{aligned}
\end{gather}

On the other hand, using $\int_{\mathbb{R}^d}|u^{TF}_\tau|^2\,dx=1$, we obtain
\begin{gather*}\label{}
	\begin{aligned} 
		\frac{1}{A_\tau^2}&=	\int_{|x|\leq\sqrt[p]{\mu^{TF}/C_0}-\tau}\eta^2_\tau(x)|u^{TF}|^2\,dx\\
		&=\int_{\mathbb{R}^d}|u^{TF}|^2\,dx+\int_{\Omega^2_\tau}\left( \eta^2_\tau(x)-1\right) |u^{TF}|^2\,dx-\int_{\Omega^1_\tau}|u^{TF}|^2\,dx\\
		&\geq1-\int_{\Omega^1_\tau\cup\Omega^2_\tau}|u^{TF}|^2\,dx\\
		&\geq1-O(\tau)\ \ \ \text{as }\tau\to0^+,
	\end{aligned}
\end{gather*}
which then implies that
\begin{gather}\label{524}
	\begin{aligned} 
		1\leq{A_\tau^2}\leq1+O(\tau)\ \ \ \text{as }\tau\to0^+.
	\end{aligned}
\end{gather}

By direct computations, we derive from \eqref{521} and \eqref{524} that
\begin{gather}\label{322}
	\begin{aligned} 
		\int_{\mathbb{R}^d}|\nabla u^{TF}_\tau|^2\,dx
		&=A_\tau^2\int_{\mathbb{R}^d}\left( |\nabla\eta_\tau u^{TF}|^2+2\nabla\eta_\tau u^{TF}\eta_\tau\nabla u^{TF}+|\eta_\tau\nabla u^{TF}|^2\right) \,dx\\
		&\lesssim\frac{1}{\tau^{2}}+\frac{1}{\tau}\int_{\Omega^2_\tau}|u^{TF}\nabla u^{TF}|\,dx+\int_{|x|\leq\sqrt[p ]{\mu^{TF}/C_0}-\tau}|\nabla u^{TF}|^2\,dx\\
  &=:\frac{1}{\tau^{2}}+\frac{1}{\tau}I_1+I_2.
	\end{aligned}
\end{gather}
Next, we give an upper bound estimate of $I_1$ and $I_2$.  Let 
\begin{gather*}
	\begin{aligned} 
		f(r):=\mu^{TF}-C_0r^p,\ \ r\in\left[0,\left(\frac{\mu^{TF}}{C_0}\right)^{\frac{1}{p}}\right].
	\end{aligned}
\end{gather*}
For any $0<r<\sqrt[p ]{\mu^{TF}/C_0}-\tau$, using the Taylor series with the Lagrange remainder, there exists a constant $\xi\in\left(r,\sqrt[p]{\mu^{TF}/C_0}\right)$ such that
\begin{gather}\label{f}
    \begin{aligned} 
    f(r)&=f\left(\left(\frac{\mu^{TF}}{C_0}\right)^{\frac{1}{p}}\right)+f'\left(\left(\frac{\mu^{TF}}{C_0}\right)^{\frac{1}{p}}\right)\left(r-\left(\frac{\mu^{TF}}{C_0}\right)^{\frac{1}{p}}\right)+\frac{f''\left(\xi\right)}{2}\left(r-\left(\frac{\mu^{TF}}{C_0}\right)^{\frac{1}{p}}\right)^2\\
    &=pC_0\left(\frac{\mu^{TF}}{C_0}\right)^{\frac{p-1}{p}}\left(\left(\frac{\mu^{TF}}{C_0}\right)^{\frac{1}{p}}-r\right)-\frac{p(p-1)C_0\xi^{p-2}}{2}\left(\left(\frac{\mu^{TF}}{C_0}\right)^{\frac{1}{p}}-r\right)^2\\
    &\geq pC_0\left(\frac{\mu^{TF}}{C_0}\right)^{\frac{p-1}{p}}\tau-2p(p-1)C_0\xi^{p-2}\tau^2\\
    &\geq O(\tau)\ \ \ \text{as }\tau\to0^+.
    \end{aligned}
\end{gather}
Moreover, we have
\begin{gather}\label{323}
	\begin{aligned} 
		I_1
		&=\frac{pC_0^2}{4}\int_{\Omega^2_\tau}\left(\mu^{TF}-C_0|x|^p\right)^{-\frac{1}{2}}|x|^{p-1}\,dx\\
		&=\frac{pC_0^2}{4}\int_{\sqrt[p]{\mu^{TF}/C_0}-2\tau}^{\sqrt[p]{\mu^{TF}/C_0}-\tau}\,dr\int_{\partial B(r)}(\mu^{TF}-C_0r^p)^{-\frac{1}{2}}r^{p-1}\,dS\\
		&=\frac{p\omega_dC_0^2}{4}\int_{\sqrt[p]{\mu^{TF}/C_0}-2\tau}^{\sqrt[p]{\mu^{TF}/C_0}-\tau}(\mu^{TF}-C_0r^p)^{-\frac{1}{2}}r^{d+p-2}\,dr\\
		&\leq C\tau^{\frac{1}{2}}\ \ \ \text{as }\tau\to0^+.
	\end{aligned}
\end{gather}
Similarly, we obtain
\begin{gather}\label{324}
	\begin{aligned} 
		I_2
		&=\frac{p^2C_0^2}{4}\int_{|x|\leq\sqrt[p]{\mu^{TF}/C_0}-\tau}\left(\mu^{TF}-C_0|x|^p\right)^{-\frac{3}{2}}|x|^{2(p-1)}\,dx\\
		&=\frac{p^2C_0^2}{4}\int_0^{\sqrt[p]{\mu^{TF}/C_0}-\tau}\,dr\int_{\partial B(0,r)}(\mu^{TF}-C_0r^p)^{-\frac{3}{2}}r^{2(p-1)}\,dS\\
		&=\frac{p^2\omega_dC_0^2}{4}\int_0^{\sqrt[p]{\mu^{TF}/C_0}-\tau}(\mu^{TF}-C_0r^p)^{-\frac{3}{2}}r^{d+2p-3}\,dr\\&
		\leq C\tau^{-\frac{3}{2}}\ \ \ \text{as }\tau\to0^+.
	\end{aligned}
\end{gather}
Thus, we deduce from \eqref{322}, \eqref{323} and \eqref{324} that 
\begin{gather}\label{hnable3}
	\begin{aligned} 
		\int_{\mathbb{R}^d}|\nabla u^{TF}_\tau|^2\,dx
		\lesssim\tau^{-2}\ \ \ \text{as }\tau\to0^+.
	\end{aligned}
\end{gather}

Now, we claim that
\begin{gather}\label{228}
	\begin{aligned} 
		0\leq E^{TF}(u^{TF}_\tau)-E^{TF}(u^{TF})\leq O(\tau)\ \ \ \text{as }\tau\to0^+.
	\end{aligned}
\end{gather}
Indeed, we deduce from \eqref{521} and \eqref{524} that
\begin{gather*}
	\begin{aligned} 
		0&\leq E^{TF}(u^{TF}_\tau)-E^{TF}(u^{TF})\\
		&=\frac{1}{2}\int_{\mathbb{R}^d}|x|^p\left(|u^{TF}_\tau|^2-|u^{TF}|^2\right)\,dx+\frac{1}{6}\int_{\mathbb{R}^d}\left(|u^{TF}_\tau|^6-|u^{TF}|^6\right)\,dx\\
		&\leq C\int_{\mathbb{R}^d}\left|A_\tau^2\eta_\tau^2(x)-1\right|\,dx+C\int_{\mathbb{R}^d}\left|A_\tau^6\eta_\tau^6(x)-1\right|\,dx\\
		&\leq C\int_{\Omega_\tau^2}\,dx+O(\tau^\frac{3}{2})\\
		&\leq O(\tau)\ \ \ \text{as }\tau\to0^+,
	\end{aligned}
\end{gather*}
where $\Omega_\tau^2$ given as \eqref{Omega}.

Next, we show that
\begin{gather}\label{428}
	\begin{aligned} 
		\left| \frac{\tau^{p}}{2}\int_{\mathbb{R}^d}V(\tau^{-1}x)|w_\tau|^2\,dx-\frac{C_0}{2}\int_{\mathbb{R}^d}|x|^p|w_\tau|^2 \,dx\right| \leq O(\tau^{p-\alpha})\ \ \ \text{as }\tau\to0^+.
	\end{aligned}
\end{gather}
We deduce from \eqref{32}, \eqref{22} and \eqref{410} that
\begin{gather}\label{4251}
	\begin{aligned} 
		\int_{\mathbb{R}^d}|x|^\alpha|w_\tau|^2\,dx&\leq
		\int_{B(R)}|x|^\alpha|w_\tau|^2\,dx+\int_{\mathbb{R}^d\setminus B(R)}|x|^\alpha|w_\tau|^2\,dx\\
		&\leq R^\alpha\int_{B(R)}|w_\tau|^2\,dx+\int_{\mathbb{R}^d\setminus B(R)}\tau^{p}|\tau^{-1}x|^p|w_\tau|^2\,dx\\
		&\leq R^\alpha+{R}\int_{\mathbb{R}^d\setminus B(R)}\tau^{p}V(\tau^{-1}x)|w_\tau|^2\,dx\\
		&\lesssim1\quad \text{as }\tau\to0^+,
	\end{aligned}
\end{gather}
where using $R>1$ (see also \eqref{32}). On the other hand, using by \textbf{(V$_4$)}, we have
\begin{gather}\label{426}
	\begin{aligned} 
		\lim_{\tau\to0^+}\frac{V(\tau^{-1}x)-C_0|\tau^{-1}x|^p}{|\tau^{-1}x|^\alpha}=C_2, \qquad \text{in } \mathbb{R}^d\setminus\{0\}.
	\end{aligned}
\end{gather}
where $C_0$ and $C_1$ are given constants (see \eqref{188} and \eqref{144}).	Thus, we have
\begin{gather*}
	\begin{aligned} 
		\left| \frac{\tau^{p}}{2}\int_{\mathbb{R}^d}V(\tau^{-1}x)|w_\tau|^2\,dx-\frac{C_0}{2}\int_{\mathbb{R}^d}|x|^p|w_\tau|^2 \,dx\right| 
		&\leq\frac{C_2\tau^{p-\alpha}}{2}\int_{\mathbb{R}^d}|x|^\alpha|w_\tau|^2 \,dx\\
		&\leq O(\tau^{p-\alpha})\ \ \ \text{as }\tau\to0^+.
	\end{aligned}
\end{gather*}

Now, we prove \eqref{hclaim} holds. Combining \eqref{410}, \eqref{428} and the fact that $E^{TF}(u^{TF})\leq E^{TF}(w_\tau)$, we deduce that 
\begin{gather}\label{similarly}
	\begin{aligned} 
		E^{TF}(u^{TF})&\leq\frac{\tau^{p+2}}{2}\int_{\mathbb{R}^d}|\nabla w_\tau|^2\,dx
		+\frac{C_0}{2}\int_{\mathbb{R}^d}|x|^p|w_\tau|^2 \,dx+\frac{1}{6}\int_{\mathbb{R}^d}|w_\tau|^6\,dx\\
		&\leq I_{\tau}(w_\tau)+\frac{\tau^{\frac{p}{2}}}{4}\int_{\mathbb{R}^d}|w_\tau|^4\,dx+\frac{C_0}{2}\int_{\mathbb{R}^d}|x|^p|w_\tau|^2 \,dx-\frac{\tau^{p}}{2}\int_{\mathbb{R}^d}V(\tau^{-1}x)|w_\tau|^2\,dx\\
		&\leq I_{\tau}(w_\tau)+O\left( \tau^{p-\alpha}\right)+O\left( \tau^{\frac{p}{2}}\right) \ \ \ \text{as }\tau\to0^+.
	\end{aligned}
\end{gather}
Similar to \eqref{similarly}, we deduce from \eqref{hnable3} that
\begin{gather*}
	\begin{aligned} 
		I_{\tau}(u^{TF}_\tau)%&=\frac{\tau^{p+2}}{2}\int_{\mathbb{R}^d}|\nabla u^{TF}_\tau|^2\,dx +\frac{\tau^{p}}{2}\int_{\mathbb{R}^d}V(\tau^{-1}x)|u^{TF}_\tau|^2\,dx \pm\frac{\tau^{\frac{p}{2}}}{4}\int_{\mathbb{R}^d}|u^{TF}_\tau|^4\,dx+\frac{1}{6}\int_{\mathbb{R}^d}|u^{TF}_\tau|^6\,dx\\
		&\leq E^{TF}(u^{TF}_\tau)+\frac{\tau^{p+2}}{2}\int_{\mathbb{R}^d}|\nabla u^{TF}_\tau|^2\,dx+\frac{\tau^{\frac{p}{2}}}{4}\int_{\mathbb{R}^d}|u^{TF}_\tau|^4\,dx\\
		&\quad+\frac{\tau^{p}}{2}\int_{\mathbb{R}^d}V(\tau^{-1}x)|u^{TF}_\tau|^2\,dx-\frac{C_0}{2}\int_{\mathbb{R}^d}|x|^p|u^{TF}_\tau|^2\,dx\\
		&\leq E^{TF}(u^{TF}_\tau)+O(\tau^{p})+O(\tau^{p-\alpha})\ \ \ \text{as }\tau\to0^+.
	\end{aligned}
\end{gather*}	
Using the above two estimates and $I_{\tau}(w_\tau)\leq I_{\tau}(u^{TF})$, one can deduce  that
\begin{gather*}
	\begin{aligned} 
		E^{TF}(u^{TF})&\leq  I_{\tau}(w_\tau)+O\left( \tau^{p-\alpha}\right)+O\left( \tau^{\frac{p}{2}}\right)\\
		&\leq E^{TF}(u^{TF}_\tau)+O\left( \tau^{p-\alpha}\right)+O\left( \tau^{\frac{p}{2}}\right)\\
		&\leq E^{TF}(u^{TF})+O\left( \tau^{p-\alpha}\right)+O(\tau)\ \ \ \text{as }\tau\to0^+,
	\end{aligned}
\end{gather*}
which implies that \eqref{hclaim} holds. 

Meanwhile, if $w_\tau$ is a non-negative minimizer of $e(\tau)$, we deduce from \eqref{410} and \eqref{428}   that
\begin{gather*}
	\begin{aligned} 
		\frac{\tau^{p+2}}{2}\int_{\mathbb{R}^d}|\nabla w_\tau|^2\,dx&=I_\tau(w_\tau)-E^{TF}(w_\tau)-\frac{\kappa\tau^{\frac{p}{2}}}{4}\int_{\mathbb{R}^d}|w_\tau|^4\,dx\\
		&\quad+\frac{\tau^{p}}{2}\int_{\mathbb{R}^d}V(\tau^{-1}x)|u^{TF}_\tau|^2\,dx-\frac{C_0}{2}\int_{\mathbb{R}^d}|x|^p|u^{TF}_\tau|^2\,dx\\
		&\leq \left| I_\tau(w_\tau)-E^{TF}(u^{TF})\right| +\frac{\tau^{\frac{p}{2}}}{4}\int_{\mathbb{R}^d}|w_\tau|^4\,dx+\frac{C_2\tau^{p-\alpha}}{2}\int_{\mathbb{R}^d}|x|^\alpha|w_\tau|^2 \,dx\\
		&\leq\left| e(\tau)-e^{TF}(\infty)\right| +O(\tau^{\frac{p}{2}})+O(\tau^{p-\alpha})\\
		&\leq O(\tau^\sigma)\ \ \ \text{as }\tau\to0^+.
	\end{aligned}
\end{gather*}
We therefore obtain \eqref{4241} and the proof of the lemma is finished.
\end{pf}

\subsection{Thomas-Fermi Limit}\label{limit}%Limiting behavior for Normalized ground states}
A limit often considered in the literature, and called the Thomas-Fermi limit, occurs when the kinetic energy is small in front of the trapping and interaction terms. Note that, we see from \eqref{4241} that the blow-up rate of kinetic energy is smaller than the decay rate of its coefficients $\tau^{p+2}$ as $\tau\to0^+$.

Firstly, we prove the following lemma.

\begin{lem}\label{lemmntau}
Assume that $1\leq d\leq3$ and $V(x)$ satisfies \textbf{(V$_1$)}--\textbf{(V$_4$)}.  Let $\sigma:=\min\{p-\alpha,1\}$ we then have
\begin{gather}\label{cg3}
	\begin{aligned} 
		|\mu_\tau-\mu^{TF}|\leq O(\tau^\sigma)\ \ \ \text{as }\tau\to0^+.
	\end{aligned}
\end{gather} 
Here $\mu_\tau$ and $\mu^{TF}$ were defined as \eqref{hwmu}  and \eqref{tfmu}.
\end{lem}
\begin{pf} First, we claim that
\begin{gather}\label{cg2}
	\begin{aligned} I_\tau(w_\tau)=E^{TF}(w_\tau)+O(\tau^\sigma)\qquad\text{as}\ \ \tau\to0^+.
	\end{aligned}
\end{gather} 
Similar to \eqref{similarly}, using by \eqref{428}, we obtain
\begin{gather}\label{4281}
	\begin{aligned} 
		E^{TF}(w_\tau)&=\frac{C_0}{2}\int_{\mathbb{R}^d}|x|^p|w_\tau|^2 \,dx+\frac{1}{6}\int_{\mathbb{R}^d}|w_\tau|^6\,dx\\
		&\leq I_{\tau}(w_\tau)+\frac{\tau^{\frac{p}{2}}}{4}\int_{\mathbb{R}^d}|w_\tau|^4\,dx+\frac{C_0}{2}\int_{\mathbb{R}^d}|x|^p|w_\tau|^2 \,dx-\frac{\tau^{p}}{2}\int_{\mathbb{R}^d}V(\tau^{-1}x)|w_\tau|^2\,dx\\
		%		&\leq I_{\tau}(w_\tau)+\frac{\tau^{\frac{p}{2}}}{4}\int_{\mathbb{R}^d}|w_\tau|^4\,dx+\frac{C_2\tau^{p-\alpha}}{2}\int_{\mathbb{R}^d}|x|^\alpha|w_\tau|^2 \,dx\\
		&\leq I_{\tau}(w_\tau)+O\left( \tau^{p-\alpha}\right)+O\left( \tau^{\frac{p}{2}}\right) \ \ \ \text{as }\tau\to0^+.
	\end{aligned}
\end{gather}
On the other hand,  we deduce from  \eqref{410}, \eqref{hnable3} and  \eqref{428} that
	\begin{align*} 
		I_{\tau}(w_\tau)&\leq\frac{\tau^{p+2}}{2}\int_{\mathbb{R}^d}|\nabla u^{TF}_\tau|^2\,dx +\frac{\tau^{p}}{2}\int_{\mathbb{R}^d}V(\tau^{-1}x)|u^{TF}_\tau|^2\,dx \pm\frac{\tau^{\frac{p}{2}}}{4}\int_{\mathbb{R}^d}|u^{TF}_\tau|^4\,dx+\frac{1}{6}\int_{\mathbb{R}^d}|u^{TF}_\tau|^6\,dx\\
		%\leq I_{\tau}(u^{TF}_\tau)
		&\leq E^{TF}(u^{TF}_\tau)+\frac{\tau^{p+2}}{2}\int_{\mathbb{R}^d}|\nabla u^{TF}_\tau|^2\,dx+\frac{\tau^{\frac{p}{2}}}{4}\int_{\mathbb{R}^d}|u^{TF}_\tau|^4\,dx\\
		&\quad+\frac{\tau^{p}}{2}\int_{\mathbb{R}^d}V(\tau^{-1}x)|u^{TF}_\tau|^2\,dx-\frac{C_0}{2}\int_{\mathbb{R}^d}|x|^p|w_\tau|^2\,dx\\
		&\leq E^{TF}(u^{TF}_\tau)+\frac{\tau^{p+2}}{2}\int_{\mathbb{R}^d}|\nabla u^{TF}_\tau|^2\,dx+\frac{\tau^{\frac{p}{2}}}{4}\int_{\mathbb{R}^d}|u^{TF}_\tau|^4\,dx+\frac{C_2\tau^{p-\alpha}}{2}\int_{\mathbb{R}^d}|x|^\alpha|u^{TF}_\tau|^2 \,dx\\
		&\leq E^{TF}(u^{TF}_\tau)+O(\tau^{p})+O(\tau^{p-\alpha})+O(\tau^{\frac{p}{2}})\\
		&\leq E^{TF}(u^{TF})+O(\tau^{p-\alpha})+O(\tau^{\frac{p}{2}})+O(\tau)\\
		&\leq E^{TF}(w_\tau)+O(\tau^{p-\alpha})+O(\tau) \ \ \ \text{as }\tau\to0^+.		
	\end{align*}
Moreover, combining \eqref{4281}, we then have 
\begin{gather}
	\begin{aligned} E^{TF}(w_\tau)&\leq I_\tau(w_\tau)+O\left( \tau^{p-\alpha}\right)+O\left( \tau^{\frac{p}{2}}\right)\\
		&\leq E^{TF}(w_\tau)+O(\tau^{p-\alpha})+O(\tau) \ \ \ \text{as }\tau\to0^+,
	\end{aligned}
\end{gather}
which implies \eqref{cg2} holds.

Next, we claim that
\begin{gather}\label{435}
	\begin{aligned} 
		\left|\int_{\mathbb{R}^d} |w_\tau|^6-|u^{TF}|^6 \,dx\right|\leq O(\tau^\sigma) \ \ \ \text{as }\tau\to0^+.
	\end{aligned}
\end{gather}
Using by \textbf{(V$_3$)}, we have
\begin{gather}\label{432}
	\begin{aligned}
		&\quad\lim_{\tau\to0^+}\frac{(\nabla V(\tau^{-1}x)\cdot \tau^{-1}x)-pC_0|\tau^{-1}x|^p}{|\tau^{-1}x|^\alpha}\\
		&=\lim_{|y|\to+\infty}\frac{(\nabla V(y)\cdot y)-pC_0|y|^p}{|y|^\alpha}=C_1.
	\end{aligned}
\end{gather} 
Thus, it follows from \eqref{4251} and \eqref{432} that
\begin{gather}\label{00}
	\begin{aligned}
		&\quad\left|\frac{pC_0}{2}\int_{\mathbb{R}^d}|x|^p|w_\tau|^2\,dx-\frac{\tau^{p-1}}{2}\int_{\mathbb{R}^d}(\nabla V(\tau^{-1}x)\cdot x)|w_\tau|^2\,dx\right| \\
		&\leq O(\tau^{p-\alpha})\int_{\mathbb{R}^d}|x|^\alpha|w_\tau|^2\,dx\\
		&\leq O(\tau^{p-\alpha})\quad \text{as }\tau\to0^+.
	\end{aligned}
\end{gather} 
On the one hand, we deduce from \eqref{hIPI} that
\begin{gather}\label{249}
	\begin{aligned}
		&\quad\frac{p\tau^{p}}{2}\int_{\mathbb{R}^d}V(\tau^{-1}x)|u|^2\,dx-\frac{d}{3}\int_{\mathbb{R}^d}|w_\tau|^6\,dx\\
		&={\tau^{p+2}}\int_{\mathbb{R}^d}|\nabla w_\tau|^2\,dx +\frac{\kappa\tau^{\frac{p}{2}} d}{4}\int_{\mathbb{R}^d}|w_\tau|^4\,dx\\
		&\quad+\frac{pC_0}{2}\int_{\mathbb{R}^d}|x|^p|w_\tau|^2\,dx-\frac{\tau^{p-1}}{2}\int_{\mathbb{R}^d}(\nabla V(\tau^{-1}x)\cdot x)|w_\tau|^2\,dx\\
		&\quad+\frac{p\tau^{p}}{2}\int_{\mathbb{R}^d}V(\tau^{-1}x)|u|^2\,dx-\frac{pC_0}{2}\int_{\mathbb{R}^d}|x|^p|w_\tau|^2\,dx.
	\end{aligned}
\end{gather}

By \eqref{410},  we have
\begin{gather}\label{cg}
	\begin{aligned} 
		0\leq\tau^{\frac{p}{2}}\int_{\mathbb{R}^d}|w_\tau|^4\,dx\leq O(\tau^{\frac{p}{2}})\quad \text{as }\tau\to0^+.
	\end{aligned}
\end{gather}

Moreover, we deduce from \eqref{00}--\eqref{cg} and \eqref{428} that
\begin{gather*}%\label{00}
	\begin{aligned}
		&\quad\lim_{\tau\to0^+}\left( \frac{p\tau^{p}}{2}\int_{\mathbb{R}^d}V(\tau^{-1}x)|u|^2\,dx-\frac{d}{3}\int_{\mathbb{R}^d}|w_\tau|^6\,dx\right)\\ &=\lim_{\tau\to0^+}\left({\tau^{p+2}}\int_{\mathbb{R}^d}|\nabla w_\tau|^2\,dx +\frac{\kappa\tau^{\frac{p}{2}} d}{4}\int_{\mathbb{R}^d}|w_\tau|^4\,dx\right) \\
		&\quad+\lim_{\tau\to0^+}\left(\frac{pC_0}{2}\int_{\mathbb{R}^d}|x|^p|w_\tau|^2\,dx-\frac{\tau^{p-1}}{2}\int_{\mathbb{R}^d}(\nabla V(\tau^{-1}x)\cdot x)|w_\tau|^2\,dx\right) \\
		&\quad+\lim_{\tau\to0^+}\left(\frac{p\tau^{p}}{2}\int_{\mathbb{R}^d}V(\tau^{-1}x)|w_\tau|^2\,dx-\frac{pC_0}{2}\int_{\mathbb{R}^d}|x|^p|w_\tau|^2\,dx\right) \\
		&=0.
	\end{aligned}
\end{gather*} 
Therefore, we obtain that, passing if necessary to a subsequence,
\begin{align}\label{cg1}
	\lim_{\tau\to0^+}\frac{p\tau^{p}}{2}\int_{\mathbb{R}^d}V(\tau^{-1}x)|u|^2\,dx=\lim_{\tau\to0^+}\frac{d}{3}\int_{\mathbb{R}^d}|w_\tau|^6\,dx.
\end{align}
By \eqref{hclaim} and \eqref{cg2}, we deduce that
\begin{gather*}
	\begin{aligned} 	\lim_{\tau\to0^+}|E^{TF}(w_\tau)-E^{TF}(u^{TF})|
		\leq\lim_{\tau\to0^+}|E^{TF}(w_\tau)-I_\tau(w_\tau)|+\lim_{\tau\to0^+}|I_\tau(w_\tau)-E^{TF}(u^{TF})| =0.
	\end{aligned}
\end{gather*} 
Combining \eqref{hTFPI}, \eqref{hclaim} and \eqref{cg1} 
\begin{gather*}
	\begin{aligned} 	\lim_{\tau\to0^+}\frac{2d+p}{6p}\int_{\mathbb{R}^d}|w_\tau|^6=\lim_{\tau\to0^+}E^{TF}(w_\tau)=E^{TF}(u^{TF})
		=\frac{2d+p}{6p}\int_{\mathbb{R}^d}|u^{TF}|^6\,dx.
	\end{aligned}
\end{gather*}
This implies that
\begin{gather}\label{cg4}
	\begin{aligned} 
		\lim_{\tau\to0^+}\int_{\mathbb{R}^d}|w_\tau|^6\,dx=\int_{\mathbb{R}^d}|u^{TF}|^6\,dx.
	\end{aligned}
\end{gather}

On the other hand, since $w_\tau$ and $u^{TF}$ are minimizers of $e(\tau)$ and $e^{TF}(\infty)$, thus, by the definition of $e(\tau)$ and $e^{TF}(\infty)$, we deduce from \eqref{hIPI} and \eqref{hTFPI} that
\begin{gather*}
	\begin{aligned} 
		\frac{2d+p}{6p}\int_{\mathbb{R}^d}|u^{TF}|^6\,dx=e^{TF}(\infty)
	\end{aligned}
\end{gather*}
and
\begin{gather}
	\begin{aligned}
		\frac{2d+p}{6p}\int_{\mathbb{R}^d}|w_\tau|^6\,dx&=e(\tau)-\frac{(p+2)\tau^{p+2}}{2p}\int_{\mathbb{R}^d}|\nabla w_\tau|^2\,dx-\frac{\kappa(p+d)\tau^{\frac{p}{2}}}{4p}\int_{\mathbb{R}^d}|w_\tau|^4\,dx  \\ &\quad	-\left( \frac{\tau^{p}}{2}\int_{\mathbb{R}^d}V(\tau^{-1}x)|w_\tau|^2\,dx-\frac{\tau^{p-1}}{2p}\int_{\mathbb{R}^d}(\nabla V(\tau^{-1}x)\cdot x)|w_\tau|^2\,dx\right).
	\end{aligned}
\end{gather}
Thus, combining \eqref{410}, \eqref{hclaim}, \eqref{4241}, \eqref{428} and \eqref{00}, we have
\begin{gather*}
	\begin{aligned} 
		&\quad\frac{2d+p}{6p}\left|\int_{\mathbb{R}^d} |w_\tau|^6-|u^{TF}|^6 \,dx\right|\\
		&\leq\left|e(\tau)-e^{TF}(1)\right| +\frac{(p+2)\tau^{p+2}}{2p}\int_{\mathbb{R}^d}|\nabla w_\tau|^2\,dx+\frac{(p+d)\tau^{\frac{p}{2}}}{4p}\int_{\mathbb{R}^d}|w_\tau|^4\,dx\\
		&\quad+\left| \frac{\tau^{p}}{2}\int_{\mathbb{R}^d}V(\tau^{-1}x)|w_\tau|^2\,dx-\frac{pC_0}{2}\int_{\mathbb{R}^d}|x|^p|w_\tau|^2\,dx\right|\\
		&\quad+\left| \frac{pC_0}{2}\int_{\mathbb{R}^d}|x|^p|w_\tau|^2\,dx-\frac{\tau^{p-1}}{2p}\int_{\mathbb{R}^d}(\nabla V(\tau^{-1}x)\cdot x)|w_\tau|^2\,dx \right|\\
		&\leq O(\tau^\sigma)\qquad \text{as }\tau\to0^+.
	\end{aligned}
\end{gather*}
This implies \eqref{435} holds. Moreover, we deduce from \eqref{hwmu}  and \eqref{tfmu} that
\begin{gather*}
	\begin{aligned} 
		\left|\mu_\tau-\mu^{TF}\right|&\leq2\left|e(\tau)-e^{TF}(\infty)\right| +\frac{2}{3}\left|  \int_{\mathbb{R}^d}|w_\tau|^6-|u^{TF}|^6\,dx\right|	+	\frac{\tau^{\frac{p}{2}}}{2}\int_{\mathbb{R}^d}|w_\tau|^4\,dx\\
		&\leq O(\tau^\sigma)\qquad \text{as }\tau\to0^+.
	\end{aligned}
\end{gather*}
Thus, \eqref{cg3} holds.
\end{pf}

\noindent
\emph{\textbf{Proof of the Theorem \ref{the2}:}}
First, we prove that \eqref{h1} holds. Using by \eqref{410}, $\{w_\tau\}$ is uniformly bounded in $L^2(\mathbb{R}^d)\cap L^6(\mathbb{R}^d)$ as $\tau\to0^+$ and there exists a
subsequence $\{w_{\tau_k}\}$ such that
\begin{gather}\label{624}
    \begin{aligned}
        w_{\tau_k}\rightharpoonup v\ \ \ \text{weakly in } L^2(\mathbb{R}^d)\cap L^6(\mathbb{R}^d).
    \end{aligned}
\end{gather}
We \textbf{claim that}
\begin{gather}\label{not0}
    \begin{aligned} 
        \int_{\mathbb{R}^d}|v|^2\,dx\not=0.	
    \end{aligned}
\end{gather} 
Indeed, by elliptic regularity and  bootstrapping the argument, we obtain that  $\varphi_{N}(x)\in C(\mathbb{R}^d)$, and $\varphi_{N}(x)$ is a radially symmetric decreasing positive function. Thus, $0$ is the unique maximum point of $\varphi_{N}(x)$ and also the maximum point of $w_\tau(x)$ by \eqref{hscaling}. Thus, we deduce from  \eqref{410} that 
\begin{gather*}
    \begin{aligned} 
	1\lesssim\int_{\mathbb{R}^d}|w_\tau|^6\,dx\leq w^4_\tau(0)\int_{\mathbb{R}^d}|w_\tau|^2\,dx=w^4_\tau(0).
    \end{aligned}
\end{gather*}
This implies that  
\begin{gather}\label{wtau1}
    \begin{aligned} 
	w_\tau(0)\gtrsim 1.
    \end{aligned}
\end{gather} 
Since $w_\tau\in C(\mathbb{R}^d)$ (the same regularity as $\varphi_{N}$), it follows that \eqref{not0} holds.

Note that $\Delta w_\tau(0)<0$. In view of \eqref{hELw}, \eqref{mutf} and  applying Lemma \ref{lemmntau}, for $\tau$ small enough, maximun $w_\tau(0)>0$ satisfies the following inequality
\begin{gather*}
    \begin{aligned} 
	w_\tau^5(0)
	\lesssim  w_\tau(0)+ \tau w_\tau^3(0).
    \end{aligned}
\end{gather*}
Hence, we have
\begin{gather}\label{wbound}
    \begin{aligned} 
	w_\tau(0)
	\lesssim1.
    \end{aligned}
\end{gather}
Since $0$ is the maximum point of $w_\tau(x)$, combining \eqref{wtau1} with \eqref{wbound}, we have
\begin{gather*}
    \begin{aligned} 
	\|\varphi_{N}\|_{L^{\infty}}\sim N^{-\frac{d}{2d+p}}.
    \end{aligned}
\end{gather*}	
Thus, \eqref{624} holds and this implies that \eqref{h3}.

In view of  \eqref{cg3} and \eqref{624}, taking the weakly limit to  \eqref{hELw}, we deduce that
\begin{gather*}
\begin{aligned} 
	v\mu^{TF}=C_0|x|^pv+v^5.
\end{aligned}
\end{gather*}
Thus, by \eqref{utf} and \eqref{not0}, we obtain
\begin{gather*}
\begin{aligned} 
	v=\left[\mu^{TF}-C_0|x|^p\right]_+^{\frac{1}{4}}=u^{TF}.
\end{aligned}
\end{gather*}
This means that $\|v\|^2_{2}=\|u^{TF}\|^2_{2}=1$, combining \eqref{cg4}, we then deduce that
\begin{gather*}
\begin{aligned} 
	w_{\tau_k}\to u^{TF} \ \ \text{strongly in}\ \ L^2(\mathbb{R}^d)\cap L^6(\mathbb{R}^d).
\end{aligned}
\end{gather*}
It is obvious to that
\begin{gather*}
\begin{aligned} 
	w_{\tau_k}\to u^{TF} \ \ \text{strongly in}\ \ L^q(\mathbb{R}^d)\ \ \text{for }2\leq q\leq 6.
\end{aligned}
\end{gather*}
Thus, \eqref{h1} holds.

Finally, we claim that \eqref{h2}. Combining \eqref{h0310}, \eqref{hEtau0} and \eqref{hclaim}, we have
\begin{gather*}
\begin{aligned} 
	\lim_{N\to+\infty}\frac{\mathscr{E}(N)}{N^{\frac{2p}{2d+p}}}=\lim_{\tau\to0^+}I_\tau(w_\tau)=E^{TF}(u^{TF})=\frac{C_0}{2}\int_{\mathbb{R}^d}|x|^p|u^{TF}|^2\,dx+\frac{1}{6}\int_{\mathbb{R}^d}|u^{TF}|^6\,dx.
\end{aligned}
\end{gather*}
Direct calculation, we obtain
\begin{gather}\label{}
\begin{aligned} 
	\int_{\mathbb{R}^d}|x|^p|u^{TF}|^2\,dx=&\int_{|x|\leq\sqrt[p]{\mu^{TF}/C_0}}|x|^p(\mu^{TF}-C_0|x|^p)^{\frac{1}{2}}\,dx\\
	=&\int_{0}^{\sqrt[p]{\mu^{TF}/C_0}}\,dr\int_{\partial B(r)}r^p(\mu^{TF}-C_0r^p)^{\frac{1}{2}}\,dS\\
	%		=&\omega_d\int_{0}^{\sqrt[p]{\mu^{TF}/C_0}}(\mu^{TF}-C_0r^p)^{\frac{1}{2}}r^{d+p-1}\,dr\\
	%		=&\frac{\omega_d}{p}\int_{0}^{\mu^{TF}/C_0}(\mu^{TF}-C_0t)^{\frac{1}{2}}t^{\frac{d}{p}}\,dt\\
	=&\frac{\omega_dC_0^{-\frac{d+p}{p}}}{p}\int_{0}^{\mu^{TF}}(\mu^{TF}-t)^{\frac{1}{2}}t^{\frac{d}{p}}\,dt\\
	=&\frac{\omega_dC_0^{-\frac{d+p}{p}}\left(\mu^{TF}\right)^{\frac{2d+3p}{2p}}}{p}\int_{0}^{\mu^{TF}}(1-\frac{t}{\mu^{TF}})^{\frac{1}{2}}\left(\frac{t}{\mu^{TF}}\right)^{\frac{d}{p}}\,d\left(\frac{t}{\mu^{TF}} \right)\\
	=&\frac{\omega_dC_0^{-\frac{d+p}{p}}\left(\mu^{TF}\right)^{\frac{2d+3p}{2p}}}{p}\mathcal{B}\left(\frac{d+p}{p},\frac{3}{2} \right),
\end{aligned}
\end{gather}
Combining \eqref{hTFPI}, we obtain
\begin{gather*}
\begin{aligned} 
	\int_{\mathbb{R}^d}|u^{TF}|^6\,dx=\frac{3pC_0}{2d}	\int_{\mathbb{R}^d}|x|^2|u^{TF}|^2\,dx=\frac{3\omega_dC_0^{-\frac{d}{p}}\left(\mu^{TF}\right)^{\frac{2d+3p}{2p}}}{2d}\mathcal{B}\left(\frac{d+p}{p},\frac{3}{2} \right).
\end{aligned}
\end{gather*}
Therefore, we obtain
\begin{gather*}
\begin{aligned} 
	\lim_{N\to+\infty}\frac{\mathscr{E}(N)}{N^{\frac{2p}{2d+p}}}=\frac{\omega_d(2d+p)C_0^{-\frac{d}{p}}\left(\mu^{TF}\right)^{\frac{2d+3p}{2p}}}{4pd}\mathcal{B}\left( \frac{d+p}{p},\frac{3}{2}\right).
\end{aligned}
\end{gather*}
Then \eqref{h2} holds. This completes the proof of Theorem \ref{the2}.  \hfill$\Box$

\section{$L^\infty$-Convergence Rate of Ground States}

The following results are interpolation estimates in the spirit of the Gagliardo-Nirenberg inequalities (see e.g. \cite{Nirenbergelliptic} or \cite[Lemma A.1]{CV1993Asymptotics})
\begin{pro}\label{nablees}
Assume $u$ satisfies
\begin{align*}
	-\Delta u=f\quad\text{in}\ K\subset\mathbb{R}^d.
\end{align*}
Then
\begin{align*}
	|\nabla u(x)|^2\leq
	C\left(\|f\|_{L^\infty(K)}\|u\|_{L^\infty(K)}+\frac{1}{dist^2(x,\partial K)}\|u\|_{L^\infty(K)}^2 \right) 
	\quad\text{in}\ K\subset\mathbb{R}^d,
\end{align*}
where $C$ is a constant depending only on $d$.
\end{pro}

Based on the above proposition, we are ready to prove the estimation of $\|\Delta w_\tau\|_{L^{\infty}}$.
\begin{lem}\label{lem44}
For any  $x\in B\left( \sqrt[p]{\mu^{TF}/C_0}\right)$, we have
\begin{align}\label{60}
	|\Delta w_\tau(x)|\leq O(\tau^{-\frac{3p+6}{4}})\quad\text{as}\ \tau\to0^+.
\end{align}
\end{lem}
\begin{pf}
Firstly, we deduce from \eqref{hELw} that
\begin{align}\label{61}
	-\Delta  w_\tau(x)=f(x)\quad\text{in}\ \mathbb{R}^d,
\end{align}
where 
$$f_\tau(x):={\tau^{-(p+2)}}\left[ -\tau^{p}V(\tau^{-1}x)w_\tau(x)-\kappa\tau^{\frac{p}{2}} w_\tau^3(x)-w_\tau^5(x)+\mu_\tau w_\tau(x)\right]. $$
Combining Lemma \ref{lemmntau}  and \eqref{wbound}, we deduce that
\begin{gather}\label{43}
	\begin{aligned} 
		\|w_\tau\|_{L^\infty(\mathbb{R}^d)}\lesssim1\quad\text{as}\ \tau\to0^+
	\end{aligned}
\end{gather}
and 
\begin{gather*}
	\begin{aligned} 
		\|f_\tau\|_{L^\infty\left(B(\sqrt[p]{\mu^{TF}/C_0}+2)\right)}\lesssim\tau^{-(p+2)}\quad\text{as}\ \tau\to0^+.
	\end{aligned}
\end{gather*}
Hence, it follow from Proposition \ref{nablees} that for any $x\in B \left( \sqrt[p]{\mu^{TF}/C_0}+2\right)$,
\begin{align}\label{453}
	|\nabla w_\tau(x)|^2\leq C\left(\frac{1}{\tau^{p+2}}+\frac{1}{dist^2\left( x,\partial B \left( \sqrt[p]{\mu^{TF}/C_0}+2\right)\right) } \right),
\end{align}
where $C>0$ is a constant depending only on $d$. Therefore, we taking $\tau>0$ small enough, such that for any $x\in B\left(  \sqrt[p]{\mu^{TF}/C_0}+1\right)$, we have $\tau^{p+2}\leq dist^2(x,\partial B \left( \sqrt[p]{\mu^{TF}/C_0}+2\right))$. Thus, we deduce from \eqref{453} that for any $x\in  B\left(  \sqrt[p]{\mu^{TF}/C_0}+1\right)$,
\begin{align}\label{63}
	|\nabla w_\tau(x)|\lesssim\tau^{-\frac{p+2}{2}}\quad\text{as}\ \tau\to0^+.
\end{align}

For any $x\in  B \left( \sqrt[p]{\mu^{TF}/C_0}+1\right)$, by \eqref{61}, we get
\begin{align*}
	-\Delta(\nabla w_\tau)& = -{\tau^{-3}}\nabla V(\tau^{-1}x)\cdot w_\tau+\tau^{-(p+2)}\left[ -\tau^{p}V(\tau^{-1}x)-3\kappa\tau^{\frac{p}{2}} w_\tau^2-5w_\tau^4+\mu_\tau\right]\nabla w_\tau\\
    &=:I_\tau^1+I_\tau^2.
\end{align*}
Next, we estimate $I_\tau^1+I_\tau^2$ in two cases: $x=0$ and $x\in  B \left( \sqrt[p]{\mu^{TF}/C_0}+1\right)\setminus\{0\}$.

\emph{Case 1: $x=0$.}   It follows easily from $V\in C^1$ that there exists a constant $C$ such that for any $\tau>0$,
\begin{align*}
\nabla V(\tau^{-1}x)\leq C\quad \text{and}\quad V(\tau^{-1}x)\leq C.
\end{align*}
Thus,
\begin{align*}
	\|I_\tau^1+I_\tau^2\|_{L^\infty\left( \{0\}\right) }\lesssim\tau^{-\frac{3p+6}{2}}\quad\text{as}\ \tau\to0^+.
\end{align*}

\emph{Case 2: $x\in  B \left( \sqrt[p]{\mu^{TF}/C_0}+1\right)\setminus\{0\}$.} From \textbf{(V$_1$)} and \textbf{(V$_4$)} that 
\begin{align*}
	\lim_{\tau\to0^+}\frac{V(\tau^{-1}x)-C_0|\tau^{-1}x|^p}{|\tau^{-1}x|^\alpha}=\lim_{|r|\to\infty}\frac{V(r)-C_0|r|^p}{|r|^\alpha}=C_2,
\end{align*}
where using $V$ is radially symmetric. This implies that
\begin{align}\label{46}
V(\tau^{-1}x)=C_0|\tau^{-1}x|^p+C_2|\tau^{-1}x|^\alpha+o(1)\quad\text{as}\ \tau\to0^+.
\end{align}
Combining \eqref{43}, \eqref{63} and \eqref{46}, we obtain
\begin{align*}
	\|I_\tau^2\|_{L^\infty\left( B \left( \sqrt[p]{\mu^{TF}/C_0}+1\right)\setminus\{0\}\right) }\lesssim\tau^{-\frac{3p+6}{2}}\quad\text{as}\ \tau\to0^+.
\end{align*}
Similarly, we deduce from \textbf{(V$_3$)} that for any $x\in B\left( \sqrt[p]{\mu^{TF}/C_0}+1\right)\setminus\{0\}$,
\begin{align}\label{47}
\tau^{-1}\nabla V(\tau^{-1}x)\cdot \tau^{-1}x=C_0p|\tau^{-1}x|^p+C_1|\tau^{-1}x|^\alpha+o(1)\quad\text{as}\ \tau\to0^+.
\end{align}
Combining \eqref{43}, \eqref{63} and \eqref{47}, we obtain
\begin{align*}
	\|I_\tau^1\|_{L^\infty\left( B \left( \sqrt[p]{\mu^{TF}/C_0}+1\right)\setminus\{0\}\right) }\lesssim\tau^{-(p+1)}\quad\text{as}\ \tau\to0^+.
\end{align*}
Thus,
\begin{align*}
	\|I_\tau^1+I_\tau^2\|_{L^\infty\left( B \left( \sqrt[p]{\mu^{TF}/C_0}+1\right)\setminus\{0\}\right) }\lesssim\tau^{-\frac{3p+6}{2}}\quad\text{as}\ \tau\to0^+.
\end{align*}
To sum up, we have
\begin{align*}
	\|I_\tau^1+I_\tau^2\|_{L^\infty\left( B \left( \sqrt[p]{\mu^{TF}/C_0}+1\right)\right) }\lesssim\tau^{-\frac{3p+6}{2}}\quad\text{as}\ \tau\to0^+.
\end{align*}

Again using Proposition \ref{nablees} we get
\begin{align*}
	|\nabla (\nabla w_\tau)|^2&\lesssim\left(\frac{1}{\tau^{\frac{3p+6}{2}}}+\frac{1}{dist^2\left(x,\partial B  \left( \sqrt[p]{\mu^{TF}/C_0}+1\right)\right)} \right),\quad x\in B \left( \sqrt[p]{\mu^{TF}/C_0}+1\right).
\end{align*}
Hence, for any $x\in B \left( \sqrt[p]{\mu^{TF}/C_0}\right)$, we get
\begin{align*}
	|\Delta w_\tau(x)|\lesssim\tau^{-\frac{3p+6}{4}}\quad\text{as}\ \tau\to0^+.
\end{align*}
Thus, \eqref{60} holds and the proof of the lemma is finished.
\end{pf}

For the rest of this subsection, we assume $$f(r)=\mu^{TF}-C_0r^p,\ \ r\in\left[0,\sqrt[p ]{\mu^{TF}/C_0}\right].$$  Based on the above lemma, we are ready to prove Theorem \ref{the3}.

\noindent
\emph{\textbf{Proof of the Theorem \ref{the3}:}} Note that \eqref{114} and \eqref{115} in Theorem \ref{the3} comes from Lemma \ref{lemmntau} and \ref{lem44}.

Notice that $w_\tau(x)\geq0$ and $u^{TF}(x)\geq0$ satisfy the following equation
\begin{gather}\label{71}
\begin{aligned} 
	-\tau^{p+2}\Delta  w_\tau(x)+\tau^pV(\tau^{-1}x)w_\tau(x)+\kappa\tau^{\frac{p}{2}} w_\tau^3(x)+w_\tau^5(x)
	=\mu_\tau w_\tau(x)\ \ \text{in }\mathbb{R}^d
\end{aligned}
\end{gather}
and
\begin{gather}\label{72}
\begin{aligned} 
	C_0|x|^pu^{TF}(x)+\left(u^{TF}(x)\right)^{5}=\mu^{TF}u^{TF}(x)\ \ \text{in }\mathbb{R}^d,
\end{aligned}
\end{gather}
respectively. Hence, we have
\begin{gather*}
\begin{aligned} 
	&\quad-\tau^{p+2}\Delta  w_\tau+\kappa \tau^{\frac{p}{2}} w_\tau^3+C_0|x|^p\left( w_\tau-u^{TF}\right)+\left(\tau^pV(\tau^{-1}x)-C_0|x|^p \right) w_\tau+\left[ w_\tau^5-\left(u^{TF}\right)^{5}\right]\\
	&=\mu_\tau w_\tau-\mu^{TF}u^{TF}\\
	&=\left( \mu_\tau-\mu^{TF}\right)  w_\tau+\mu^{TF}\left(w_\tau-u^{TF}\right)\quad\text{in }\mathbb{R}^d.
\end{aligned}
\end{gather*}
Moreover, we deduce from \eqref{utf} that for any $x\in B\left( \sqrt[p]{\mu^{TF}/C_0}\right) $,
\begin{gather}\label{53}
\begin{aligned} 
	&\quad-\tau^{p+2}\Delta  w_\tau+\kappa\tau^{\frac{p}{2}} w_\tau^3+\left(\tau^pV(\tau^{-1}x)-C_0|x|^p \right) w_\tau+\left[ w_\tau^5-\left(u^{TF}\right)^{5}\right]\\
	&	=\left(\mu_\tau-\mu^{TF}\right)w_\tau+\left(u^{TF}\right)^{4}\left( w_\tau-u^{TF}\right) .
\end{aligned}
\end{gather}

From the binomial theorem, we obtain
\begin{gather}\label{55}
\begin{aligned} 
	w_\tau^5-\left(u^{TF}\right)^{5}&=\left(\left( w_\tau-u^{TF}\right) +u^{TF}\right)^{5}-\left(u^{TF}\right)^{5}\\
	&=\left( w_\tau-u^{TF}\right)^5
	+5\left( w_\tau-u^{TF}\right)^4u^{TF}
	+10\left( w_\tau-u^{TF}\right)^3\left(u^{TF}\right)^{2}\\
	&\quad+10\left( w_\tau-u^{TF}\right)^2\left(u^{TF}\right)^{3}
	+5\left( w_\tau-u^{TF}\right)\left(u^{TF}\right)^{4}.
\end{aligned}
\end{gather}
Thus, for any $x\in B\left( \sqrt[p]{\mu^{TF}/C_0}\right) $, we have
\begin{gather}\label{}
\begin{aligned} 
	4\left( w_\tau-u^{TF}\right)\left(u^{TF}\right)^{4}   &=-\left( w_\tau-u^{TF}\right)^5
	-5\left( w_\tau-u^{TF}\right)^4u^{TF}
	-10\left( w_\tau-u^{TF}\right)^3\left(u^{TF}\right)^{2}\\
	&\quad-10\left( w_\tau-u^{TF}\right)^2\left(u^{TF}\right)^{3}+\tau^{p+2}\Delta  w_\tau-\kappa\tau^{\frac{p}{2}} w_\tau^3+\left[\mu_\tau-\mu^{TF}\right]w_\tau\\
	&\quad-\left(\tau^pV(\tau^{-1}x)-C_0|x|^p \right) w_\tau.
\end{aligned}
\end{gather}
Moreover, we deduce from \eqref{426}, \eqref{cg3} and \eqref{60} that for any $x\in B\left( \sqrt[p]{\mu^{TF}/C_0}\right)$,
\begin{gather}\label{56}
\begin{aligned} 
	\left(u^{TF}\right)^{4}  \left\|w_\tau-u^{TF}\right\|_{L^\infty} &\lesssim\tau^\sigma
	+\left\|w_\tau-u^{TF}\right\|^2_{L^\infty},
\end{aligned}
\end{gather}
where $\sigma=\min\{p-\alpha,1\}$.

We can prove (i) of Theorem \ref{the3}.  For any compact subset $K$ of $B(\sqrt[p]{\mu^{TF}/C_0})$, we obtain $u^{TF}(x)\sim1$ for $x\in K$. Moreover, for any $x\in K$, we have
\begin{gather}\label{581}
\begin{aligned} 
	\left(u^{TF}\right)^{8} -4C\tau^\sigma\gtrsim1-\tau^\sigma>0\quad\text{as}\ \tau\to0^+.
\end{aligned}
\end{gather}
Thus, we deduce from \eqref{56} and \eqref{581} that
\begin{gather}\label{}
\begin{aligned} 
	\left\|w_\tau-u^{TF}\right\|_{L^\infty\left( K\right) }
 &\leq\frac{\left(u^{TF}\right)^{4}-\sqrt{\left(u^{TF}\right)^{8}-4C\tau^\sigma}}{2C}\\
	&=\frac{2\tau^\sigma}{\left(u^{TF}\right)^{4}+\sqrt{\left(u^{TF}\right)^{8}-4C\tau^\sigma}}\\
	&\lesssim\tau^\sigma\quad\text{as}\ \tau\to0^+.
\end{aligned}
\end{gather}

Next, we prove (ii) of Theorem \ref{the3}. Similar to \eqref{f}, for any $|x|<\sqrt[p]{\mu^{TF}/C_0}-\left( \tau|\ln \tau|\right)^\epsilon$, there exists a constant $\xi\in(|x|,\sqrt[p]{\mu^{TF}/C_0}-\left( \tau|\ln \tau|\right)^\epsilon)$ such that
\begin{gather}\label{f2}
    \begin{aligned} 
    f(r)&=pC_0\left(\frac{\mu^{TF}}{C_0}\right)^{\frac{p-1}{p}}\left(\left(\frac{\mu^{TF}}{C_0}\right)^{\frac{1}{p}}-r\right)-\frac{p(p-1)C_0\xi^{p-2}}{2}\left(\left(\frac{\mu^{TF}}{C_0}\right)^{\frac{1}{p}}-r\right)^2\\
    &\geq pC_0\left(\frac{\mu^{TF}}{C_0}\right)^{\frac{p-1}{p}}\left( \tau|\ln \tau|\right)^\epsilon-2p(p-1)C_0\xi^{p-2}\left( \tau|\ln \tau|\right)^{2\epsilon}\\
    &\gtrsim\left( \tau|\ln \tau|\right)^\epsilon\ \ \ \text{as }\tau\to0^+,
    \end{aligned}
\end{gather}
where $\epsilon>0$ is a constant independent of $\tau$. For any $x\in B\left( \sqrt[p]{\mu^{TF}/C_0}-\left( \tau|\ln \tau|\right)^\epsilon\right)$, we deduce form \eqref{utf} and \eqref{f2} that
\begin{gather*}
\begin{aligned} 
	\left( \mu^{TF}\right)^{\frac{1}{4}}\geq u^{TF}(x)\geq\left(f(|x|)\right)^{\frac{1}{4}}\gtrsim\left( \tau|\ln \tau|\right)^\frac{\epsilon}{4}\quad\text{as}\ \tau\to0^+.
\end{aligned}
\end{gather*}
Let $0<\epsilon\leq\sigma/2$, for any $x\in B\left( \sqrt[p]{\mu^{TF}/C_0}-\left( \tau|\ln \tau|\right)^\epsilon\right)$, we have
\begin{gather}\label{58}
\begin{aligned} 
	\left(u^{TF}\right)^{8} -4C\tau^\sigma \gtrsim\left( \tau|\ln \tau|\right)^{2\epsilon}-\tau^\sigma\gtrsim\left( \tau|\ln \tau|\right)^{2\epsilon}>0\quad\text{as}\ \tau\to0^+.
\end{aligned}
\end{gather}
Thus, we deduce from \eqref{56} and \eqref{58} that for any $0<\epsilon\leq\sigma/2$,
\begin{gather}\label{}
\begin{aligned} 
	\left\|w_\tau-u^{TF}\right\|_{L^\infty\left( B\left( \sqrt[p]{\mu^{TF}/C_0}-\left( \tau|\ln \tau|\right)^\epsilon\right)\right) }
 &\leq\frac{\left(u^{TF}\right)^{4}-\sqrt{\left(u^{TF}\right)^{8}-4C\tau^\sigma}}{2C}\\
	&=\frac{2\tau^\sigma}{\left(u^{TF}\right)^{4}+\sqrt{\left(u^{TF}\right)^{8}-4C\tau^\sigma}}\\
	&\lesssim\tau^{\sigma-\epsilon}|\ln \tau|^{-\epsilon}\quad\text{as}\ \tau\to0^+.
\end{aligned}
\end{gather}

Finally, we prove \eqref{117}.
Combining  \eqref{71} and $w_\tau(x)\geq0$, for any $\tau>0$
\begin{gather*}\label{}
\begin{aligned} 
	-\Delta w_\tau(x)+\frac{\tau^pV(\tau^{-1}x)-\mu_\tau}{\tau^{p+2}}w_\tau(x)-\tau^{-\frac{4+p}{2}}w_\tau^3(x)\leq0,\quad x\in\mathbb{R}^d.
\end{aligned}
\end{gather*}
Moreover, we deduce that 
\begin{gather*}\label{}
\begin{aligned} 
	-\Delta w_\tau(x)+\frac{\tau^pV(\tau^{-1}x)-C_0|x|^p}{\tau^{p+2}}w_\tau(x)+\frac{C_0|x|^p-\mu_\tau}{\tau^{p+2}}w_\tau(x)-\tau^{-\frac{4+p}{2}}w_\tau^3(x)\leq0,\quad x\in\mathbb{R}^d.
\end{aligned}
\end{gather*}
Thus, we find that as $\tau\to0^+$
\begin{gather}\label{422}
\begin{aligned} 
	-\Delta w_\tau(x)+\frac{C_1|x|^\alpha}{\tau^{\alpha+2}}w_\tau(x)+\frac{C_0|x|^p-\mu_\tau}{2\tau^{p+2}}w_\tau(x)+\tau^{-\frac{4+p}{2}}\left( \frac{C_0|x|^p-\mu_\tau}{2\tau^{\frac{p}{2}}}-w_\tau^2(x)\right) w_\tau(x)\leq0
\end{aligned}
\end{gather}
for any $ x\in\mathbb{R}^d$.

On the one hand, for any $x\in \mathbb{R}^d\setminus B\left( \sqrt[p]{\mu^{TF}/C_0}+ \tau^{{p}/{2}-\epsilon}|\ln \tau|^{-\epsilon}\right)$, using the Taylor series with the Lagrange remainder,  there exists a constant $\xi\in\left(\sqrt[p]{\mu^{TF}/C_0}+ \tau^{{p}/{2}-\epsilon}|\ln \tau|^{-\epsilon},|x|\right)$ such that
\begin{gather*}
    \begin{aligned} 
  C_0|x|^p-\mu^{TF}&=pC_0\left({\frac{\mu^{TF}}{C_0}}\right)^{\frac{p-1}{p}}\left(|x|-\left({\frac{\mu^{TF}}{C_0}}\right)^{\frac{1}{p}}\right)-\frac{p(p-1)C_0\xi^{p-2}}{2}\left(\left({\frac{\mu^{TF}}{C_0}}\right)^{\frac{1}{p}}-|x|\right)^2\\
    &\geq pC_0\left({\frac{\mu^{TF}}{C_0}}\right)^{\frac{p-1}{p}}\tau^{\frac{p}{2}-\epsilon}|\ln \tau|^{-\epsilon}-2p(p-1)C_0\xi^{p-2}\tau^{{p}-2\epsilon}|\ln \tau|^{-2\epsilon}\\
    &\gtrsim\tau^{\frac{p}{2}-\epsilon}|\ln \tau|^{-\epsilon}\ \ \ \text{as }\tau\to0^+.
    \end{aligned}
\end{gather*}
Combining  \eqref{cg3} and $\sigma=\min\{p-\alpha,1\}$, for any $x\in \mathbb{R}^d\setminus B\left( \sqrt[p]{\mu^{TF}/C_0}+ \tau^{{p}/{2}-\epsilon}|\ln \tau|^{-\epsilon}\right)$, we have
\begin{gather}\label{423}
\begin{aligned} 
	\frac{C_0|x|^p-\mu_\tau}{2\tau^{\frac{p}{2}}}&\geq\frac{(C_0|x|^p-\mu^{TF})-|\mu^{TF}-\mu_\tau|}{2\tau^{\frac{p}{2}}}\\
	%&\gtrsim\frac{\tau^{\frac{p}{2}-\epsilon}|\ln \tau|^{-\epsilon}-\tau^\sigma}{\tau^{\frac{p}{2}}}\\
	&\gtrsim\frac{1}{\tau^\epsilon|\ln \tau|^{\epsilon}}-\tau^{\frac{p}{2}-\sigma}\quad(p\geq2)\\
	&\gtrsim \frac{1}{\tau^\epsilon|\ln \tau|^{\epsilon}} \quad\text{as}\ \tau\to0^+.
\end{aligned}
\end{gather}
Note that $\|w_\tau\|_{L^{\infty}}\sim1$, thus for any $x\in \mathbb{R}^d\setminus B\left( \sqrt[p]{\mu^{TF}/C_0}+ \tau^{{p}/{2}-\epsilon}|\ln \tau|^{-\epsilon}\right)$, 
\begin{gather}\label{424}
\begin{aligned} 
	\frac{C_0|x|^p-\mu_\tau}{2\tau^{p+2}}-w_\tau^2(x)\geq0 \quad\text{as}\ \tau\to0^+.
\end{aligned}
\end{gather}

On the other hand, for any $x\in \mathbb{R}^d\setminus B\left( \sqrt[p]{\mu^{TF}/C_0}+ \tau^{{p}/{2}-\epsilon}|\ln \tau|^{-\epsilon}\right)$, we deduced from \eqref{423} that  there exists a constant $\beta>0$ such that
\begin{gather}\label{425}
\begin{aligned} 
	\frac{C_0|x|^p-\mu_\tau}{2\tau^{\frac{p}{2}}}\geq\frac{\beta^2}{4}.
\end{aligned}
\end{gather}
Combining \eqref{422}, \eqref{424} and \eqref{425}, for any $x\in \mathbb{R}^d\setminus B\left( \sqrt[p]{\mu^{TF}/C_0}+ \tau^{{p}/{2}-\epsilon}|\ln \tau|^{-\epsilon}\right)$, we know that  $w_\tau(x)\geq0$ satisfies
\begin{gather*}
\begin{aligned} 
	-\Delta w_\tau(x)+\frac{\beta^2\tau^{-\frac{p+4}{2}}}{4}w_\tau(x)\leq0 \quad\text{as}\ \tau\to0^+.
\end{aligned}
\end{gather*}

Let $v\in L^{\infty}\left( \mathbb{R}^d\setminus B\left( \sqrt[p]{\mu^{TF}/C_0}+ \tau^{{p}/{2}-\epsilon}|\ln \tau|^{-\epsilon}\right),\mathbb{R}\right) $ be such that
\begin{gather*}
\begin{cases}
	-\Delta v(x)+\frac{\beta^2\tau^{-\frac{p+4}{2}}}{4}v(x)=0\ \ &\text{if}\ x\in\mathbb{R}^d\setminus B\left( \sqrt[p]{\mu^{TF}/C_0}+ \tau^{{p}/{2}-\epsilon}|\ln \tau|^{-\epsilon}\right),\\
	v(x)=w_{\tau}(x)&\text{if}\ x\in B\left( \sqrt[p]{\mu^{TF}/C_0}+ \tau^{{p}/{2}-\epsilon}|\ln \tau|^{-\epsilon}\right),\\
	\lim_{|x|\to\infty}v(x)=0.
\end{cases}
\end{gather*}
By \cite[Lemma 6.4]{JFAvm} with $W={\beta^2\tau^{-\frac{p+4}{2}}}/{4}$, for any $x\in\mathbb{R}^d\setminus B\left( \sqrt[p]{\mu^{TF}/C_0}+ \tau^{{p}/{2}-\epsilon}|\ln \tau|^{-\epsilon}\right)$, there exists $C\in\mathbb{R}^+$ such that 
\begin{align*}
v(x)\leq C|x|^{-\frac{d-1}{2}}e^{-\frac{\beta\tau^{-\frac{p+4}{4}}}{2}|x|}.
\end{align*}
Hence, by the comparison principle \cite{Okavian}, for any $x\in\mathbb{R}^d\setminus B\left( \sqrt[p]{\mu^{TF}/C_0}+ \tau^{{p}/{2}-\epsilon}|\ln \tau|^{-\epsilon}\right)$, we then have
\begin{align*}
w_\tau(x)\leq v(x)\lesssim|x|^{-\frac{d-1}{2}}e^{-\frac{\beta\tau^{-\frac{p+4}{4}}}{2}|x|} \quad\text{as}\ \tau\to0^+.
\end{align*}
This completes the proof of Theorem \ref{the3}. \hfill$\Box$

\appendix
\renewcommand{\appendixname}{Appendix A~\Alph{section}}

\section{Existence of Ground States}\label{secA}

The following lemma is well known, one can see its proof in  \cite[Theorem XIII.6.7]{Reed-simon-book} or \cite{TSPjz}.
\begin{lem}\label{compact}
Suppose $1\leq d\leq 3$ and $V(x)\in L^{\infty}_{loc}(\mathbb{R}^d)$ satisfies $\lim_{|x|\to\infty}V(x)=\infty$.	For $2\leq r <2^*$,  the embedding $\mathscr{H}_N\hookrightarrow L^r (\mathbb{R}^d)$ is compact, where $2^*=\infty$ for $d=1,2$, or $2^*=6$ for $d=3$.
\end{lem}

\begin{lem}[{\cite{JMPcp}}]\label{xpingf}
Assume that $V(x)$ is increasing and nonnegative radially symmetric. Let $u^*$ be the symmetric rearrangement of $u$ that vanishes at infinity, then
\begin{gather}
	\begin{aligned} 
		\int_{\mathbb{R}^d}V|u^*|^2\,dx\leq\int_{\mathbb{R}^d}V|u|^2\,dx.
	\end{aligned}
\end{gather}
\end{lem}

\noindent
\emph{\textbf{Proof of the Theorem \ref{existence}:}}
Recall from \cite{anailsis}, we have $|\nabla |u||\leq|\nabla u|$, a.e., $x\in\mathbb{R}^d$. Thus, we deduce that $E(u)\geq E(|u|)$. Therefore, it suffices to consider the real non-negative minimizers of \eqref{amin}.

For any $N > 0$, let $\{u_n\}\subset\mathscr{H}_N$ be a minimizing sequence of $e(N)$. On the one hand, if $\kappa=+1$, it is easy fined that  
\begin{gather}\label{1121}
\begin{aligned}
	E(u_n)\geq\frac{1}{2}\int_{\mathbb{R}^d}\left( |\nabla u_n|^2+V(x)|u_n|^2\right) \,dx.
\end{aligned}
\end{gather}
On the other hand, if $\kappa=-1$,  we deduce from H\"older inequality and Young inequality that
\begin{gather*}
\begin{aligned} 
	\int_{\mathbb{R}^d}|u_n|^4\,dx&\leq\left( \int_{\mathbb{R}^d}|u_n|^6\,dx\right)^{\frac{1}{2}}\left( \int_{\mathbb{R}^d}|u_n|^2\,dx\right)^{\frac{1}{2}}\\
	&\leq\frac{2}{3}\int_{\mathbb{R}^d}|u_n|^6\,dx+\frac{3N}{8}.
\end{aligned}
\end{gather*}
Moreover, we have
\begin{gather}\label{1122}
\begin{aligned}
	E(u_n)=&\frac{1}{2}\int_{\mathbb{R}^d}\left( |\nabla u_n|^2+V(x)|u_n|^2\right) \,dx-\frac{1}{4}\int_{\mathbb{R}^d}|u_n|^4\,dx+\frac{1}{6}\int_{\mathbb{R}^d}|u_n|^6\,dx\\
	\geq&\frac{1}{2}\int_{\mathbb{R}^d}\left( |\nabla u_n|^2+V(x)|u_n|^2\right) \,dx-\frac{3N}{32}.
\end{aligned}
\end{gather}
Using by \eqref{1121} and \eqref{1122}, we then have the minimizing sequence $\{u_n\}$ is uniformly bounded in $\mathscr{H}_N$ for $\kappa=\pm1$.
Therefore, up to a subsequence $u_n\rightharpoonup u_0$ weakly in $\mathscr{H}_N$. By Lemma \ref{compact} above, we obtain
\begin{gather*}
\begin{aligned}
	u_n\to u_0\ \text{ strongly in } L^r(\mathbb{R}^d)\ \text{ for }2\leq r<2^*.
\end{aligned}
\end{gather*}
Thus, we conclude that $\|u_0\|_2^2=N$ and $E(u_0)=e(N)$, by weak lower semicontinuity. This implies $u_0$ is a minimizer of $e(N)$.

Finally, we show that the non-negative ground state  $\varphi_{N}(x)$ is radially symmetric and decreasing with respect to $|x|$. Let $\varphi_{N}^*(x)$ be the symmetric-decreasing rearrangement of $\varphi_{N}(x)$. Thus, we deduce from Lemma \ref{xpingf} and \textbf{(V$_1$)} that
\begin{gather*}
\begin{aligned} 
	\int_{\mathbb{R}^d}V(x)|\varphi_{N}^*(x)|^2\,dx\leq	\int_{\mathbb{R}^d}V(x)|\varphi_{N}(x)|^2\,dx.
\end{aligned}
\end{gather*}
We know from \cite[Chapter 3 and Lemma 7.17]{anailsis} that for any $s\geq1$,
\begin{gather*}
\begin{aligned} 
	\int_{\mathbb{R}^d}|\varphi_{N}^*(x)|^s\,dx
	=\int_{\mathbb{R}^d}|\varphi_{N}(x)|^s\,dx\ \text{ and}\ \int_{\mathbb{R}^d}|\nabla \varphi_{N}^*(x)|^2\,dx
	\leq\int_{\mathbb{R}^d}|\nabla \varphi_{N}(x)|^2\,dx.
\end{aligned}
\end{gather*}	 
It follows that
\begin{gather}\label{z2}
\begin{aligned} 
	e(N)\leq E(\varphi_{N}^*)\leq E(\varphi_{N})=e(N). 	
\end{aligned}
\end{gather}
Thus $\varphi_{N}^*(x)$ is a ground state of $e(N)$, this implies that $\varphi_{N}(x)$ is radially symmetric and decreasing with respect to $|x|$.    
\hfill$\Box$

\vspace{8mm}
\noindent\textbf{Acknowledgments.} 
The authors declare that they have no conflict of interest.\vspace{3mm}

\noindent\textbf{Data Availability Statemen.} Data sharing not applicable to this article as no datasets were generatedor analyzed during the current study.

\end{document}